\documentclass[dvips]{amsart}
\usepackage{amssymb,amsxtra,latexsym,graphics,pictex,epsfig}

\newcommand{\vanish}[1]{}
\newcommand{\F}{\mathcal{F}}

\newcommand{\A}{\mathcal{A}}
\newcommand{\B}{\mathcal{B}}
\renewcommand{\L}{\mathcal{L}}

\newcommand{\C}{\mathcal{C}}

\newcommand{\R}{\mathbb{R}}

\newcommand{\Z}{\mathbb{Z}}
\newcommand{\fq}{\mathbb{F}_q}
\newcommand{\abs}[1]{\lvert#1\rvert}

\newcommand{\supp}{\operatorname{supp}}

\newcommand{\dimm}{\operatorname{dim}}

\newcommand{\type}{\operatorname{type}}
\newcommand{\End}{\operatorname{End}}

\newcommand{\onto}{\twoheadrightarrow}
\newcommand{\into}{\hookrightarrow}
\newcommand{\iso}{\cong}
\newcommand{\isoto}{\stackrel{\cong}{\longrightarrow}}

\theoremstyle{plain}

\newtheorem{proposition}{Proposition}

\theoremstyle{definition}

\theoremstyle{remark}
\newtheorem*{remark}{Remark}

\newcommand{\aigner}{aigner97:_combin}

\newcommand{\riffle}{bayerdiaconis92:_trail}
\newcommand{\bid}{bidigare97:_hyper}

\newcommand{\bhr}{bidigarehanlon98}
\newcommand{\bbd}{billera99:_random}

\newcommand{\red}{bjoerner93:_om}

\newcommand{\brown}{brown89:_build}
\newcommand{\ken}{brown00:_semig_markov}
\newcommand{\bd}{browndiaconis:_random}
\newcommand{\borovik}{borovik95:_coxet}

\newcommand{\drep}{diaconis88:_group}

\newcommand{\dfp}{diaconisfill92:_top_to_random}

\newcommand{\gb}{grove85:_finit}

\newcommand{\humphreys}{humphreys90:_reflec_coxet}

\newcommand{\ot}{orlikterao92:_book}

\newcommand{\schar}{scharlau95:_build}

\newcommand{\solomon}{solomon76:_mackey}

\newcommand{\stanley}{stanley97:_enumer1}

\newcommand{\tits}{tits74:_build_bn}

\newcommand{\wachswhite}{wachs91:_stirl}

\newcommand{\zie}{ziegler95:_lectures}

\begin{document}
\title{Shuffles on Coxeter groups}
\author{Swapneel Mahajan}

\address{Department of Mathematics\\
Cornell University\\
Ithaca, NY 14853}
\email{swapneel@math.cornell.edu}


\begin{abstract}
The random-to-top and the riffle shuffle are two well-studied methods
for shuffling a deck of cards. These correspond to the symmetric group
$S_n$, i.e., the Coxeter group of type $A_{n-1}$. In this paper, we
give analogous shuffles for the Coxeter groups of type $B_n$ and
$D_n$. These can be interpreted as shuffles on a ``signed'' deck 
of cards. With these examples as motivation, we abstract the notion of a
shuffle algebra which captures the connection between the algebraic
structure of the shuffles and the geometry of the Coxeter groups.
We also briefly discuss the generalisation to buildings which leads to 
$q$-analogues.
\end{abstract}

\maketitle

\section{Introduction}

In a recent work, Ken Brown~\cite{\ken} used algebraic methods to
analyse random walks on a class of semigroups called ``left-regular
bands''. These walks include the hyperplane
chamber walks of Bidigare, Hanlon, and Rockmore~\cite{\bhr}. 
In this paper, we look at the special case of reflection arrangements
that arise in the study of Coxeter groups.
The random walks that we look at can be thought of as shuffles on the
Coxeter group.
The motivating examples are the riffle shuffle and the random-to-top
shuffle on a deck of $n$ cards.	
These correspond to the symmetric group
$S_n$, i.e., the Coxeter group of type $A_{n-1}$. 

The only Coxeter groups that we deal with are the ones of type $A_{n-1}$,
$B_n$ and $D_n$. 
To make this paper accessible to readers unfamiliar with Coxeter
groups we have included two appendices. 
Appendix~\ref{app:hyperplane} reviews the facts we need about hyperplane
arrangements. In Appendix~\ref{app:Coxeter}, we give a brief review of Coxeter
groups and then explain everything in
concrete terms for the three cases mentioned above.
For the general theory of Coxeter groups, we refer the reader to
\cite{\brown,\gb,\humphreys,\tits}. 

\subsection{The random walk and the method of analysis} \label{subs:rw}
Let $\Sigma$ be the simplicial (Coxeter) complex associated to a
Coxeter group $W$. Also let $\C$ be the set of chambers (or maximal
simplices) of $\Sigma$. Then $\Sigma$ is a semigroup containing $\C$
as an ideal. The product in $\Sigma$ is given by the projection maps
as explained in Appendix~\ref{app:hyperplane}.
We now describe the walk.
Let $\{w_x\}_{x\in \Sigma}$ be a probability distribution
on $\Sigma$. 
If the walk is in chamber $c \in \C$, then move to the chamber $xc$,
where $x \in \Sigma$ is chosen with probability $w_x$.
The product $xc$ is again a chamber because $\C$ is an ideal in $\Sigma$.

Next we describe the walk more algebraically.
Fix a commutative ring $k$ and consider the semigroup algebra
$k\Sigma$. The $k$-module $k\C$ spanned by the chambers
is an ideal in $k\Sigma$.  In particular, it is a module over 
$k\Sigma$; we therefore obtain a homomorphism
\[
k\Sigma\to \End_k(k\C),
\]
the latter being the ring of $k$-endomorphisms of $k\C$.  This map is
in fact an inclusion; so we can regard elements of $k\Sigma$ as
operators on $k\C$.  
Let $\sigma = \sum_{x\in \Sigma} w_x x$ be such that 
$\sum_{x \in \Sigma} w_x = 1$.
It is straightforward to check that the transition matrix of the random
walk determined by~$\{w_x\}$ is simply the matrix of the operator
``left multiplication by~$\sigma$''.  
In other words, it is the image of $\sigma$ under the inclusion
$k\Sigma \into \End_k(k\C)$.

A way to analyse the random walk is to focus on $k\Sigma$ and
understand the structure
of the subalgebra $k[\sigma]$ that $\sigma$ generates in $k\Sigma$.
Because of the inclusion map above, we can also analyse the structure
of $k[\sigma]$ in $\End_k(k\C)$.
We will make use of both viewpoints in our examples.
The reader may be more familiar with this algebraic approach in the context of
groups~\cite{\drep} rather than semigroups~\cite{\ken}.

\subsection{Nature of our examples} \label{subs:ne}
In our examples, the element $\sigma$ will never be normalised. 
In other words, if we write $\sigma = \sum_{x\in \Sigma} w_x x$
then $\sum_{x \in \Sigma} w_x \not= 1$.
So, strictly speaking, we will not have a probability distribution on $\Sigma$.
The main reason for doing this is to simplify the algebra.
So everytime we describe the random walk associated to $\sigma$,
the description will be correct upto a normalisation factor.
This will also be true for all other objects that will show up
in the analysis.
To get rid of this anamoly,
we adopt the following convention.
If the random walk description says
``Do xxx at random.''
then the normalisation factor is
``the number of ways of doing xxx''.
In general, the random walk description will involve
a sequence of independent random acts in which case
the normalisation factor will be the product of the individual factors.
In most examples, 
we will prefer to first define all objects of interest 
as elements of $k\Sigma$
and later to provide motivation by interpreting the associated random walks.
Making this translation will always involve the normalisation factor
defined above.

Next we mention that in our examples, the element $\sigma$ will always
lie in the subalgebra of $W$-invariants of $k\Sigma$.
We now discuss this subalgebra.
The $W$-invariants of $k\Sigma$ under the natural $W$-action form
a $k$-algebra $(k\Sigma)^W$.  As a $k$-module, $(k\Sigma)^W$ is free with
one basis element for each $W$-orbit in~$\Sigma$, that basis element
being the sum of the simplices in the orbit.  Since orbits correspond
to types of simplices, we get a basis vector
\[
\sigma_J= \sum_{F\in\Sigma_J} F
\]
for each $J\subseteq I$, where
$\Sigma_J$ is the set of simplices of
type~$J$ and $I$ is
the set of all labels or types of vertices of $\Sigma$.
Also define $\Sigma_j$ to be the set of simplices of rank $j$.

As already mentioned, in our examples $\sigma \in (k\Sigma)^W$ and
hence $k[\sigma]$ is a certain commutative subalgebra 
of $(k\Sigma)^W$. Since the element $\sigma$ is always motivated by
some shuffle considerations, we will call $k[\sigma]$ a shuffle algebra.  
Bidigare \cite{\bid}
proved that $(k\Sigma)^W$ is anti-isomorphic to Solomon's descent
algebra~\cite{\solomon}, which is a certain subalgebra of the group
algebra $kW$.  So 
shuffle algebras are anti-subalgebras of Solomon's descent
algebra. 

\subsection{Organization of the paper}
In Section~\ref{s:mt} we define shuffle algebras and provide a 
semi-simplicity criterion for their analysis.
In the next three sections,
we give examples of shuffle algebras of type
$A_{n-1}, B_n$ and $D_n$ respectively.
These sections should be read in conjunction with
Appendices~\ref{subs:a}, \ref{subs:b} and~\ref{subs:d} respectively.
Our examples are motivated by shuffle considerations 
on the one hand and the geometry of the Coxeter complex on the other.
The striking similarity among the examples can be traced to maps among the 
three Coxeter complexes.
These maps which give a unified framework for our examples 
are explained in Section~\ref{s:maps}.
This section gives us a good overall picture
of how things fit together and 
can play an important role in the further development of the theory.
We see an example of this in Section~\ref{s:q}, 
where we briefly discuss the generalisation of the 
random walks to buildings.

\section{The main tool} \label{s:mt}

This section is somewhat technical and
it is better to read it lightly at present.
You may want to go over the details later 
with a concrete example in mind.
We first abstract the notion of a shuffle algebra
and then develop the main tool for the analysis of these algebras.

\subsection{Shuffle algebras} \label{subs:defn}

We first need a preliminary definition.
An additive (resp. multiplicative) \emph{shuffle
semigroup} is a subsemigroup of the semigroup of non-negative
(resp. positive)
integers under addition (resp. multiplication).
Note that there is only one additive shuffle semigroup upto isomorphism.
However in the multiplicative case, there could be many.
Now we are ready to define a shuffle algebra.

We say that a subalgebra $A\subseteq k\Sigma$ is an additive (resp. multiplicative)  \emph{shuffle algebra} if
it satisfies the following conditions:

\begin{enumerate}
\item[(1)]
$\dimm_kA \leq \dimm \Sigma + 2$. 
\item[(2)] 
$A$ has a basis of the form $\sigma_0 =
1$,$\sigma_1$,$\sigma_2$,\ldots where $\sigma_j$ is an element of
$k\Sigma_j$, the span of simplices of rank $j$. 
\item[(3)] 
$A = k[\sigma_1]\subseteq k\Sigma$.
\item[(4)] 
$A$ contains an additive (resp. multiplicative) shuffle
semigroup $S$ as a spanning set.
\end{enumerate}
Note that the definition of a shuffle algebra $A$ depends on the
Coxeter complex $\Sigma$ under consideration. In fact, conditions
$(1),(2)$ and $(3)$ say that $A$ is in tune with the geometry of $\Sigma$.
Hence it is more correct to write a shuffle algebra as a pair
$(A,\Sigma)$.
But we will not bother with this since in our examples
the Coxeter complex $\Sigma$
will always be clear from context.

An alternate way to express condition $(4)$ is to say that
there exists an injective semigroup map $S \into A$ 
such that the image spans $A$.
In all our examples, the semigroup $S$ will have an interpretation 
in terms of card shuffles. This explains our terminology.

\begin{remark}
The above definition is mainly motivated by the infinite families
$A_{n-1}, B_n$ and $D_n$.
The significance of this definition for the sporadic Coxeter groups
of type $E_6$, $E_7$, etc. is not clear.
It is also natural to ask whether shuffle algebras can be classified.
We do not know the answer.
We would also like to point out that we have imposed the strictest
possible conditions that our examples led us to.
As a result, some of the algebras we consider are not shuffle
algebras;
for example, the shuffle double constructions in 
Sections~\ref{subs:brs} and~\ref{subs:drs}.
Hence a more flexible definition might work better,
say for the classification problem. 

\end{remark}

\subsection{A semisimplicity criterion}

Ken Brown~\cite{\ken} shows that if $k$ has characteristic 0 and $\sigma$ is a
non-negative integral linear combination in the canonical basis of
$k\Sigma$ consisting of all simplices
then $k[\sigma]$ is \emph{split-semisimple}, that is, $k[\sigma]
\iso k^n$ for $n=\dimm_kA$. To guarantee this result for an arbitrary
field he introduces an additional condition on $\sigma_1$.  In all our
examples, we assume that $k$ has characteristic 0. Also $\sigma_1$
will always be a non-negative integral linear combination in the
canonical basis of $k\Sigma$. Hence it follows by condition $(3)$ that
the shuffle algebras we consider are split-semisimple. We will not
rely on Brown's result however.  Instead we will use a semisimplicity
criterion which we now state.

\begin{proposition} \label{p:ss}
Let $A$ be an $n$-dim algebra containing a spanning semigroup $S$.
Suppose that there exist characters $\chi_1$,\ldots,$\chi_n$ of $S$
and elements $e_1$,\ldots,$e_n$ of $A$ such that for every $s \in S$,
we have
\begin{equation} \tag{*}
s = \sum_{i=1}^n \chi_i(s) e_i.
\end{equation}
Then $A \isoto k^n$ with $e_i$ as the primitive idempotents.
Also the $\chi_i$'s extend to characters of $A$.
\end{proposition}
\begin{proof}
The elements $e_1$,\ldots,$e_n$ of $A$ span $S$ which in turn
spans $A$. 
Since $\dimm_kA = n$ this implies that $e_1$,\ldots,$e_n$ form a basis for $A$. 
This yields a vector space isomorphism $\Phi : A \isoto k^n$, where
$\sum_{i=1}^n a_i e_i \mapsto (a_1,\ldots,a_n)$.
Hence for $s \in S$,
$s \mapsto (\chi_1(s),\ldots,\chi_n(s))$.
Since $\chi_1,\ldots,\chi_n$ are characters of $S$, we have
$\Phi(s_1 s_2) = \Phi(s_1) \Phi(s_2)$ for $s_1,s_2 \in S$.
The fact that $\Phi$ respects the algebra structures on a spanning set
$S$ of $A$ implies that it is in fact an algebra homomorphism.
\end{proof}

As mentioned before, this proposition will apply to all our
examples. The algebra $A$ will be a shuffle algebra and $S$ the
shuffle semigroup contained in $A$. 
The characters $\chi$ of $S$ that we will use are
$\chi(a) = c^a$ if $S$ is additive and 
$\chi(a) = a^c$ if $S$ is multiplicative. 
Here $c$ is a fixed non-negative integer.
The specific values of $c$ that we need to choose 
depend on the example at hand.
The advantage of Proposition~\ref{p:ss}
is that it is tailor-made for shuffle algebras
and it gives an explicit isomorphism of $A$ with $k^n$.

\begin{remark}
The Coxeter complexes $\Sigma$ of type $A_{n-1}, B_n$ and $D_n$ can be
described using the language of ``ordered partitions''. 
Also, chambers of $\Sigma$ can be described as a ``deck of cards''. 
This is explained in Appendices~\ref{subs:a}, \ref{subs:b} and~\ref{subs:d} 
respectively.
We follow the convention that ``left to right'' for an ordered
partition is ``top to bottom'' for a deck of cards.
It is essential to understand this language before reading the examples.
The side shuffle and riffle shuffle of type $A_{n-1}$ considered in
Sections~\ref{subs:ass} and~\ref{subs:ars} are the motivating
examples of this paper. It is also a good idea to first understand these two
basic examples before proceeding to the rest.
\end{remark}

\section{Examples of type $A_{n-1}$}

We consider three examples; the side shuffle, the two sided shuffle
and the riffle shuffle.
The first two are additive in nature while the third is
multiplicative.
The analysis for the first and third example parallels that of
Diaconis, Fill and Pitman~\cite{\dfp}
and Bayer and Diaconis~\cite{\riffle} 
respectively.
The difference is that we prefer to work with $(k\Sigma)^W$,
which is the more geometric description of Solomon's descent algebra.
Elements of the label set $I$ will be written as $s_1, s_2, \ldots, s_{n-1}$;
see Figure 2.1 in Appendix~\ref{subs:cd}.
Recall that $\sigma_J$, the sum of simplices of type~$J$, for $J \subseteq I$
form a basis for $(k\Sigma)^W$. Hence these elements will play a key
role in the analysis.

\subsection{The side shuffle}\label{subs:ass}
This is more commonly known as the random-to-top shuffle or the
Tsetlin library.
The element of interest is $\sigma = \sigma_{s_1}$, that is, it is the
sum of all vertices of type $s_1$.
In terms of ordered partitions, $\sigma_{s_1}$ is the sum of all
ordered two block partitions of $[n]$ such that the first block is a
singleton.
As explained in Section~\ref{subs:rw}, 
there is a random walk on a deck of $n$ cards
associated to $\sigma$.
It consists of removing a
card at random and replacing it on top. 
When $n=4$, for example, 
$(\{2\}, \{1,3,4\})$ is a typical summand of $\sigma$. Its product with
the deck $(\{1\}$,$\{2\}$,$\{3\}$,$\{4\})$ gives
$(\{2\}$,$\{1\}$,$\{3\}$,$\{4\})$, that is, the overall effect is to
remove the card labelled $2$ and put it on top.
Note that the normalisation factor in this case is $n$,
consistent with the convention of Section~\ref{subs:ne}.

In addition to $\sigma = \sigma_1$,
we define $\sigma_j = \sigma_{J_j}$ where $J_j = \{s_1,\ldots,s_j\}
\subseteq I$ for $j = 1,\ldots,n-1$. 
To explain in words,
$\sigma_1$ is the sum of all the vertices of type $s_1$,
$\sigma_2$ is the sum of all the edges of type $s_1 s_2$ and so on till
$\sigma_{n-1}$ which is the sum of all chambers of $\Sigma$.
As ordered partitions, $\sigma_{j}$ is the sum of all
ordered $(j+1)$ block partitions of $[n]$ such that the first $j$
blocks are singletons.
Also for convenience, we define $\sigma_n = \sigma_{n-1}$. 
Just like $\sigma_1$, 
we can describe 
the random walk on a deck of $n$ cards
associated to $\sigma_j$.

$\sigma_j :$ Choose $j$ distinct cards at random and put them on top in a
random order. 

\noindent Note that $\sigma_n$ was previously defined
in an artificial way.
But the random walk description for $\sigma_n$ makes perfect sense
and we now see the motivation in setting $\sigma_n = \sigma_{n-1}$. 
The two are identical as random walks and by our convention
have the same normalisation factor of $n!$.
From now on, we will leave out the discussion on normalisation.
We now present two methods to
analyse $k[\sigma_1]$.

\medskip
\noindent
{\bf The shuffle method}.

Put $\sigma_0 = 1$.
Define $a$-shuffles $S_a$ for $a \geq 0$ by 
\begin{equation} \label{e:shuffle}
S_a = \sum_{j=0}^n S(a,j) \sigma_j,
\end{equation}
where $S(a,j)$ are the Stirling numbers of the second kind.
They count the number of ways in which $a$ elements can be divided
into $j$ non-empty subsets. 
We make the convention that 
$S(a,j) = 0$ for $a < j$. Also
$S(0,0) = 1$ and $S(a,0) = 0$ for $a > 0$.
Note that
$S_0 = 1$ and $S_1 = \sigma_1$.
The motivation behind this definition will soon become clear.
Unlike for $\sigma_j$,
the description of $S_a$ in terms of ordered partitions is not
clear. However, the associated random walk on a deck of $n$ cards 
can be readily described.

$S_a :$ Pick a card at random and mark it. Repeat this process $a$
times. There is no restriction on the number of times a given card may get
marked. Move all the marked cards to the top in a random order.

To see the equivalence of the two definitions, note
that $S_a$ can be expressed as a sum indexed by the number
of distinct cards that get marked. 
Suppose that $j$ distinct cards get marked. By our marking scheme for
$S_a$, there are exactly $S(a,j)$ ways in which the same set of $j$
cards gets marked. This leads to equation~\eqref{e:shuffle}.
With the probabilistic interpretation for $S_a$, one can also convince
oneself that $S_a S_b = S_{a+b}$.
This relation was the main motivation behind the definition of $S_a$.
Thus we get an additive shuffle semigroup $S = \{S_a : a \geq 0 \}$
contained in $A = k[\sigma_1] = k[S_1]$.
The additivity says that the semigroup $S$ spans $A$.
We can deduce from equation~\eqref{e:shuffle} that
$\sigma_0=1$,$\sigma_1$,\ldots,$\sigma_{n-1}$ lie in $A$. 
Also they span $S$ and hence $A$.
Since we know that they are linearly independent in $(k\Sigma)^W$, 
it now follows that they form a basis for $A$. 
Hence $A$ is a shuffle algebra.

The analysis so far shows that the element $\sigma_1$ determines what 
$\sigma_2,\ldots,\sigma_{n-1}$ should be, assuming that we want
$\sigma_j \in \Sigma_j$. We will see this more
directly in the second method. Now we show that $A \iso k^n$.
We do this by moulding equation~\eqref{e:shuffle} in the shape of
equation~(*) of Proposition~\ref{p:ss}. 
The key step is to use the following explicit formula for the Stirling
numbers $S(a,j)$; see~\cite[pg 34]{\stanley}.
\begin{equation*} 
S(a,j) = \sum_{i=0}^j
(-1)^{j-i}\binom{j}{i}\frac{i^a}{j!},
\end{equation*}
with the convention that $0^0=1$.
Substituting the above expression in the formula for $S_a$ and
rearranging terms, we get for $a \geq 0$
\begin{equation*} 
S_a=\sum_{i=0}^n i^a e_i, \quad \text{where} \quad
e_i=\sum_{j=i}^n(-1)^{j-i}\binom{j}{i}\frac{\sigma_j}{j!}\,. 
\end{equation*}
Observe that $e_{n-1} = 0$ because $\sigma_{n-1}=\sigma_{n}$.
The remaining $e_i$'s are clearly non-zero.
Now apply Proposition~\ref{p:ss} to $S$ and $A$ along with the $n$
characters of $S$ given by $\chi_i(S_a) = i^a$  for $i=0,1,2,\ldots,n-2,n$.
This gives
$A \isoto k^{n}$. 
The isomorphism maps $S_a$ to $(0^a,1^a,2^a,\ldots,(n-2)^a,n^a)$.
These give the possible set of eigenvalues of $S_a$ considered as an
operator on any $A$-module.  
The formulas for the $e_i$'s are identical to those obtained in~\cite{\ken}.
The minimal polynomial for $S_a$ is given by
$x(x-1)(x-2^a)\cdots(x-(n-2)^a)(x-n^a)$.
If we expand this polynomial, substitute $x=S_a$ and use
equation~\eqref{e:shuffle} then we get some identities involving
the Stirling numbers.
\vanish
{
Now set
$S' = \{S_a : a > 0 \}$ and
$A'$ to be the $k$-span $\sigma_1$,\ldots,$\sigma_{n-1}$.
The proposition applied to $S'$ and $A'$ along with the $n-1$
characters $\chi_i(S_a) = i^a$  for $i=1,2,\ldots,n-2,n$ of $S'$ gives
$A' \isoto k^{n-1}$. 

}

\medskip
\noindent
{\bf The direct method}.

We first claim that
$\sigma_j \sigma_1 = j \sigma_j + \sigma_{j+1}$ if $j < n-1$ and 
$\sigma_{n-1} \sigma_1 = n \sigma_{n-1}$.
We see this from the description of the $\sigma$'s as ordered
partitions as follows. A typical summand of $\sigma_j$ is
$(\{1\}$,$\{2\}$,$\cdots$,$\{j\}$,$\{j+1 \cdots n\})$.
This term appears $j$ times in the product $\sigma_j \sigma_1$;
the $j$ summands of $\sigma_1$ that contribute
being the ones in which the element in the singleton block 
is one of $1,2, \ldots,j$.
This yields $j \sigma_j$. The remaining terms yield $\sigma_{j+1}$ and
the claim follows. For $j=1$, the claim says that $\sigma_1^2 = \sigma_1 +
\sigma_{2}$. Thus $\sigma_2 = \sigma_1^2 - \sigma_1$
is determined by $\sigma_1$ and the same is true for $\sigma_3$ and so on. 
More precisely, by induction we get
\begin{equation}\label{e:cob}
\sigma_{j+1} = \sigma_1(\sigma_1-1)\ldots(\sigma_1-j).
\end{equation}
This along with
$\sigma_{n-1} \sigma_1 = n \sigma_{n-1}$ implies that both
$\sigma_0=1$,$\sigma_1$,\ldots,$\sigma_{n-1}$ and
$\sigma_1^0$,$\sigma_1^1$,\ldots,$\sigma_1^{n-1}$
form a basis for $A$.
The basic relation satisfied by $\sigma_1$ is 
$\sigma_{n-1}(\sigma_1-n) = 0$, that is,
$\sigma_1(\sigma_1-1)\ldots(\sigma_1-n+2)(\sigma_1-n) = 0$.
Hence
$A \isoto \frac{k[x]}{x(x-1)\ldots(x-n+2)(x-n)}$,
where the map sends $\sigma_1$ to $x$.
This shows that $A$ is split semisimple.
Note that this time we did not use Proposition~\ref{p:ss}
to arrive at this conclusion.

We now give algebraic motivation for equation~\eqref{e:shuffle}
by rederiving it here.
Right now $\sigma_j$ is defined only for $0 \leq
j \leq n-1$.
But we may extend the definition of $\sigma_j$ to any $j$
using equation~\eqref{e:cob}. 
With this extension, one may check that 
$\sigma_n = \sigma_{n-1}$ and $\sigma_j = 0$ for
$j > n$.
We now find a formula for $\sigma_1^j$ by inverting
equation~\eqref{e:cob} formally.
Write 
\begin{equation*}
\sigma_1^a = \sum_{j=0}^a S(a,j) \sigma_j
\end{equation*}
for some constants $S(a,j)$.
Multiply both sides on the right by $\sigma_1$ and use the relation
$\sigma_j \sigma_1 = j \sigma_j + \sigma_{j+1}$ to obtain the
recursion:
$S(a,j) = j S(a-1,j) + S(a-1,j-1)$ with $S(a,a) = S(a,1) = 1$.
The recursion and the initial conditions show that 
$S(a,j)$ are the Stirling numbers of the second kind.
We now define the $a$-shuffle $S_a$ to be $\sigma_1^a$.
This gives equation~\eqref{e:shuffle}.
By our definition, it follows directly that $S_a S_b = S_{a+b}$.
Note that this time we did not rely on any 
probabilistic interpretation of $S_a$
to derive this additive relation.

\begin{remark}
Consider the map $k[x] \onto A$ which sends $x$ to $\sigma_1$.
Along with $\{x^j\}_{j \geq 0}$,
the other sequence which played a prominent role in our analysis was
$\{x_{(j)}\}_{j \geq 0}$, where
$x_{(j+1)} = x(x-1)\ldots(x-j)$.
Both these sequences are polynomial sequences of binomial type~\cite{\aigner}. 
This kind of structure seems to be common to all the additive examples.
The first sequence is always the same but the second sequence varies.
And the relation between the two gives us various analogues of 
the Stirling numbers.
It would be worthwhile to understand this better.

It is possible to write down an explicit formula for $\sigma_i
\sigma_j$. We state it here for completeness.
$$\sigma_i\sigma_j=\sum_{k=0}^{min(i,j)}k!\binom{i}{k}\binom{j}{k}
\sigma_{i+j-k}.$$
There is also a $q$-analogue of the side shuffle which we will
explain in Section~\ref{s:q}.
\end{remark}

\subsection{The two sided shuffle} \label{subs:atss}
This example is similar to the side shuffle and it will be useful to
keep the analogy in mind.
The element of interest is $\sigma = \sigma_{s_1} + \sigma_{s_{n-1}}$. 
In terms of ordered partitions, $\sigma$ is the sum of all
ordered two block partitions of $[n]$ such that either the first block
or the second block is a singleton.  
The associated random walk on a deck of cards consists of removing a
card at random and replacing it either on top or at the bottom; whence
the name. 

In addition to $\sigma = \sigma_1$,
we define $\sigma_j = \sum_{k=0}^j \binom{j}{k}\sigma_{J_{j,k}}$, where
$J_{j,k} = \{s_1,\ldots,s_k\} \cup \{s_{n-(j-k)},\ldots,s_{n-1}\}
\subseteq I$ for $j = 1,\ldots,n-1$. 
Note that $J_{j,k}$ is the union of the first $k$ and last $j-k$
elements of the label set. Hence the cardinality of $J_{j,k}$ is always $j$.
Also as ordered partitions, $\sigma_{J_{j,k}}$ is the sum of all
ordered $(j+1)$ block partitions of $[n]$ such that the first $k$ and
the last $j-k$ blocks are singletons.
As is evident from the formula, expressing $\sigma_j$ in this language
is not pleasant. However the probability description will be simple as
we will see.
Also put $\sigma_0 = 1$ and
$\sigma_n = 2 \sigma_{n-1}$.
Now define $a$-shuffles $S_a$ for $a \geq 0$ by 
\begin{equation} \label{e:shuffle2}
S_a = \sum_{j=0}^n S(a,j) \sigma_j,
\end{equation}
where $S(a,j)$ counts the number of ways in which a set of $a$ elements can be divided
into $j$ non-empty subsets, where in each subset the elements are further
divided into two subsets. 
We will call the $S(a,j)$'s the \emph{signed} Stirling numbers.
They satisfy the recursion:
$S(a,j) = 2j S(a-1,j) + S(a-1,j-1)$ with $S(a,a) = 1$ and $S(a,1) =
2^{a-1}$. An explicit formula is as under:
\begin{equation*} 
S(a,j) = \sum_{i=0}^j
(-1)^{j-i}\binom{j}{i}\frac{(2i)^a}{2^j j!}\,.
\end{equation*}
Note that
$S_0 = 1$ and $S_1 = \sigma_1$. 
These definitions are motivated by the following interpretations 
of the random walks associated to $\sigma_j$ and $S_a$.

$\sigma_j :$ Choose $j$ cards at random. Split these into two distinct
piles such that every split is equally likely. Move the cards in the
first pile to the top in a random order and those in the second
pile to the bottom in a random order. 

$S_a :$ Pick a card at random and mark it $T$ or $B$. Repeat this process $a$
times. If at any time we choose a card that has already been marked
then we overwrite that mark. Now move all cards marked $T$ to the top
and those marked $B$ to the bottom in a random order.

First check that $\sigma_n = 2 \sigma_{n-1}$.
Next note that $\sigma_j$ can be expressed as a sum depending on the
number of cards $k$ that are moved to the top.
There are $\binom{j}{k}$ ways for this to occur.
This leads us to the equation for $\sigma_j$ in terms of the
$\sigma_{J_{j,k}}$'s that we wrote earlier.
The interpretation of $S_a$ is slightly subtle. 
As for the side shuffle, express $S_a$ as a sum indexed by the number
of distinct cards that get marked. 
Suppose that $j$ distinct cards get marked. By our marking scheme for
$S_a$, there are exactly $2 S(a,j)$ ways in which the same set of $j$
cards gets marked. This is because we use distinct labels $T$
and $B$. However, this extra factor is compensated in the next step
where we use the last label on a marked card to decide whether it
goes to the top or the bottom. 

Now observe that $S_a S_b = S_{a+b}$.
This gives us an additive shuffle semigroup $S = \{S_a : a \geq 0 \}$
contained in $A = k[\sigma_1] = k[S_1]$.
We repeat the analysis for the side-shuffle to show that 
$A \isoto k^{n}$, where $S_a$ maps to $(0^a,2^a,4^a,\ldots,(2n-4)^a,(2n)^a)$. 
The formula for the primitive idempotents is given by
\begin{equation*} 
e_i = \sum_{j=i}^n(-1)^{j-i}\binom{j}{i}\frac{\sigma_j}{2^j j!}\,.
\end{equation*}
This time the relation $\sigma_{n}= 2 \sigma_{n-1}$ forces $e_{n-1} = 0$.
Also equation~\eqref{e:shuffle2} can be rewritten as
\[
S_a=\sum_{i=0}^n (2i)^a e_i.
\]
As for the side-shuffle, the direct method also works. Using the
description of the $\sigma$'s as ordered partitions, 
it is easy to see that
$\sigma_j \sigma_1 = 2j \sigma_j + \sigma_{j+1}$ if $j < n-1$ and 
$\sigma_{n-1} \sigma_1 = 2n \sigma_{n-1}$.
By induction we get
\begin{equation}\label{e:cob2}
\sigma_{j+1} = \sigma_1(\sigma_1-2)\ldots(\sigma_1-2j).
\end{equation}
This along with
$\sigma_{n-1} \sigma_1 = 2n \sigma_{n-1}$ implies that both
$\sigma_0=1$,$\sigma_1$,\ldots,$\sigma_{n-1}$ and
$\sigma_1^0=1$,$\sigma_1^1$,\ldots,$\sigma_1^{n-1}$
form a basis for $A$.
The basic relation satisfied by $\sigma_1$ is 
$\sigma_1(\sigma_1-2)\ldots(\sigma_1-2n+4)(\sigma_1-2n) = 0$.
Hence
$A \isoto \frac{k[x]}{x(x-2)\ldots(x-2n+4)(x-2n)}$.
This shows that $A$ is split semisimple.
Inverting equation~\eqref{e:cob2} formally and using the relation
$\sigma_j \sigma_1 = 2j \sigma_j + \sigma_{j+1}$ 
leads to the recursion of the signed Stirling numbers written earlier
and hence to equation~\eqref{e:shuffle2}.

\subsection{Riffle shuffle} \label{subs:ars}
A {\it riffle shuffle} is a common method people use for shuffling a
deck of cards. The Gilbert-Shannon-Reeds mathematical model for riffle
shuffling consists of first cutting the deck into two parts, making
the cut according to a binomial probability distribution, and then
interleaving the two parts in such a way that every interleaving is
equally likely. A natural generalisation is the $a$-shuffle in which
we cut the deck into $a$ parts (rather than just 2) and then interleave
them.

From our point of view, it is more natural to consider
the inverse riffle shuffle $S_2$ and more generally the inverse
$a$-shuffle $S_a$. It has the following description as a random walk on
a deck of $n$ cards. 

$S_a$ : Label each card randomly with an integer from $1$ to
$a$. Move all the cards labelled $1$ to the bottom of the deck,
preserving their relative order. Next move all the cards labelled $2$
above these again preserving their relative order and so on.

We can understand the inverse $a$-shuffle $S_a$ as an element of
$(k\Sigma)^W$. We begin with the element $\sigma_1 =
\sigma_{s_1}+\cdots+\sigma_{s_{n-1}}$; that is,
$\sigma_1$ is the sum of all the vertices in $\Sigma$.  In terms of
ordered partitions, it is the sum of all ordered two block partitions
of $[n]$.  This is closely related to the inverse riffle shuffle, in
fact, $S_2 = \sigma_1 + 2 \sigma_0$, where $\sigma_0 = 1$ is the one
block partition. We will show
this in more generality.

Let $\sigma_i$ be the sum of all the simplices of rank $i$ in $\Sigma$
for $i = 1,\ldots,n-1$. Then
define shuffles $S_a$ for $a \geq 1$ by
\begin{equation} \label{e:ars}
S_a = \sum_{j=1}^n \binom{a}{j} \sigma_{j-1}.
\end{equation}
In terms of ordered partitions, $\sigma_{j-1}$ is the sum of all
ordered $j$ block partitions
of $[n]$. To express $S_a$ in this language, we need the notion of a
weakly ordered 
partition of $[n]$. It is an ordered partition of $[n]$ in which the
blocks are allowed to be empty. Note that there are exactly 
$\binom{a}{j}$ weakly ordered $a$ block partitions that give the same
ordered $j$ block partition of $[n]$. 
It follows that $S_a$ is the sum of all weakly
ordered $a$ block partitions
of $[n]$. 
If we assign labels $1,2,\ldots,a$ to these weakly
ordered $a$ blocks such that the leftmost block is labelled $a$ and
the rightmost is labelled $1$, then
we get the random walk description for $S_a$ that we started with. 
It is clear from both viewpoints that 
$S_a S_b = S_{ab}$. 
This gives us a multiplicative shuffle semigroup $S = \{S_a : a \geq 1 \}$
contained in $A = k[\sigma_1] = k[S_2]$.
It now follows that 
$\sigma_0=1$,$\sigma_1$,\ldots,$\sigma_{n-1}$ is a basis for $A$ 
and that $A$ is a shuffle algebra.

Now we show that $A \iso k^n$.
Expand each $\binom{a}{j}$ as a polynomial in $a$ (with no constant
term):
\begin{equation*} 
\binom {a}{j} = \frac{a(a-1)\ldots(a-j+1)}{j!} =
\frac{1}{j!}\sum_{i=1}^j c_{ij} a^i\,
\end{equation*}
for some constants $c_{ij}$.
To be precise, $c_{ij}$ is the coefficient of $x^i$ in the polynomial
$x(x-1)\cdots(x-j+1)$.
Substituting the above expression in the formula for $S_a$ and
rearranging terms, we get for $a \geq 1$
\begin{equation*} 
S_a=\sum_{i=1}^n a^i e_i, \quad \text{where} \quad e_i=\sum_{j=i}^n
c_{ij}\frac{\sigma_{j-1}}{j!}\,. 
\end{equation*}
All the $e_i$'s are clearly non-zero.
Now apply Proposition~\ref{p:ss} to $S$ and $A$ along with the $n$
characters of $S$ given by $\chi_i(S_a) = a^i$  for $i=1,2,\ldots,n$.
This gives
$A \isoto k^{n}$. 
The isomorphism maps $S_a$ to $(a^1,a^2,\ldots,a^n)$.
These give the possible set of eigenvalues of $S_a$ considered as an
operator on any $A$-module.  
They are distinct for $a \neq 1$ and hence $S_a$ generates $A$ for all
$a \geq 2$. Also the minimal polynomial for $S_a$ on any faithful
$A$-module is given by $(x-a)(x-a^2)\cdots(x-a^n)$.
Write it as 
$\sum_{i=0}^n P_{n,i} x^i$.
Now if we substitute $x=S_a$ and use equation~\eqref{e:ars} then we
get identities involving binomial coefficients, namely
$\sum_{i=0}^n P_{n,i} \binom{a^i}{j} = 0$ for $j=1,2,\dots,n$.
It would be nice to have an explicit formula for the coefficients
$P_{n,i}$ also. 

The analysis that we gave for this example 
was the analogue of the shuffle method
of the previous two examples.
The direct computational method of those examples is not feasible here.
Writing a formula for $\sigma_j \sigma_1$ is not easy because
of the multiplicative nature of the example.

\section{Examples of type $B_n$}

We consider two examples, the side shuffle and the riffle shuffle.
They should not be regarded as mere generalisations of the ones
discussed so far.
In fact, we will see in Section~\ref{s:maps} that 
they are in a sense more fundamental and 
have a right to exist on their own.
In this section, the term partition always refers to a partition of
type $B_n$; see Appendix~\ref{subs:b}. 
Also elements of the label set $I$ will be written as $s_1, s_2, \ldots, s_{n}=t$;
see Figure 2.2 in Appendix~\ref{subs:cd}.

\subsection{The side shuffle}
The element of interest is $\sigma = \sigma_{s_1}$. 
In terms of partitions, $\sigma_{s_1}$ is the sum of all
three block partitions such that the first (and hence the last)
block is a 
singleton.  
The associated random walk on a deck of $n$ signed cards 
consists of removing a
card at random and replacing it on top with either the same sign or
its reverse both with equal probability. 
Note that the normalisation factor is $2n$.

In addition to $\sigma = \sigma_1$,
we define $\sigma_j = \sigma_{J_j}$ where $J_j = \{s_1,\ldots,s_j\}
\subseteq I$ for $j = 1,\ldots,n$. 
To explain in words,
$\sigma_1$ is the sum of all the vertices of type $s_1$,
$\sigma_2$ is the sum of all the edges of type $s_1 s_2$ and so on till
$\sigma_{n}$ which is the sum of all chambers of $\Sigma$.
As partitions, $\sigma_{j}$ is the sum of all
$(2j+1)$ block partitions such that the first (and hence the last) $j$
blocks are singletons.
Define $a$-shuffles for $a \geq 0$ by 
\begin{equation*} 
S_a = \sum_{j=0}^n S(a,j) \sigma_j,
\end{equation*}
where $S(a,j)$ are the signed Stirling numbers defined in
Section~\ref{subs:atss}.
The associated random walks on a deck of $n$ signed cards 
are as follows:

$\sigma_j :$ Choose $j$ cards at random. Move them to the top in a
random order. Now with equal probability either flip or do not flip
the sign of the chosen cards.  

$S_a :$ Pick a card at random and with equal probability either flip
or do not flip its sign. Repeat this process $a$ times. There is no
restriction on the number of times the same card may get picked. Move
all the picked cards to the top in a random order.

The interpretation of $S_a$ involves the same kind of subtlety as the
two sided shuffle of type $A_{n-1}$ discussed in Section~\ref{subs:atss}.
We observe that $S_a S_b = S_{a+b}$. After this,
the analysis works exactly like the two sided shuffle for type
$A_{n-1}$. The only difference is that $e_{n-1} \neq 0$ since we do
not have any relation like $\sigma_{n} = 2 \sigma_{n-1}$.
This gives $A \isoto k^{n+1}$, where $S_a$ maps to $(0^a,2^a,4^a,\ldots,(2n-2)^a,(2n)^a)$. 

Like the additive examples of type $A_{n-1}$, 
we can also use the direct method here.
Using the partition description of the $\sigma$'s,
it is easy to see that
$\sigma_j \sigma_1 = 2j \sigma_j + \sigma_{j+1}$ if $j \leq n-1$ and 
$\sigma_{n} \sigma_1 = 2n \sigma_{n}$ and the rest is similar.

\subsection{The riffle shuffle}\label{subs:brs}
This example is motivated by the riffle shuffle of type $A_{n-1}$.  In
this case, the shuffles $S_a$ split into two cases depending on the
parity of $a$. They can be described using the inverse $a$-shuffles
of Section~\ref{subs:ars} and a special shuffle $S_2$.
We will call $S_2$ the inverse {\it
signed} riffle shuffle or the inverse riffle shuffle of type $B_n$.

$S_2$: For every card, we either flip or do not flip its sign with equal
probability. The cards with unchanged signs move to the top 
in the same relative order and the
rest move to the bottom in the reverse relative order. 

$S_{2a}$: We do a usual inverse $a$-shuffle with labels
$1,2,\ldots,a$. Then within each of the $a$ blocks with a fixed label 
we do an inverse signed riffle shuffle $S_2$.

$S_{2a+1}$: We do a usual inverse $(a+1)$-shuffle with labels
$0,1,2,\ldots,a$. Then we do an inverse signed riffle shuffle $S_2$ 
on each block except the one labelled 0.

These shuffles are multiplicative. 
One can check directly that $S_{a} S_{b} = S_{ab}$ for any $a,b \geq
1$ irrespective of parity. However, we will derive this relation by
considering the shuffles $S_a$ as elements of
$(k\Sigma)^W$. It is natural to first analyse the even and odd shuffles
separately and then to put them together later.

\medskip
\noindent
{\bf The even part}.

The element of interest is $\sigma_1 = \sigma_{t}$. 
In terms of partitions, $\sigma_{t}$ is the sum of all
three block partitions such that the second (zero)
block is  empty.

If we try to analyse $\sigma_1^2$, $\sigma_1^3$ and so on then we
observe that they involve only those partitions that have an
empty zero block. In geometric language, they involve only those faces
whose type contains the letter $t$. This motivates the definition of
$\sigma_j$ which we now give. Put $\sigma_j = \sum_{\abs{J}=j,t \in J}
\sigma_J$ for $j=1,2,\ldots,n$.
In terms of partitions, $\sigma_{j}$ is the sum of all
$(2j+1)$ block partitions such that the zero
block is  empty.
In the spirit of the riffle shuffle of type $A_{n-1}$, we define shuffles
$S_{2a}$ for $a \geq 1$ by
\begin{equation}\label{e:ershuffle}
S_{2a} = \sum_{j=1}^n \binom{a}{j} \sigma_{j}.
\end{equation}
Note that $S_2 = \sigma_1$.
Also put $S_1 = \sigma_0 = 1$.
To express $S_{2a}$ in words, we use the notion of a weak 
partition that we defined in Appendix~\ref{subs:b}. It is a partition
in which the 
signed blocks are also allowed to be empty. 
Since the partitions are anti-symmetric, there are exactly 
$\binom{a}{j}$ weak $(2a+1)$ block partitions that give the same
$(2j+1)$ block partition. 
It follows that $S_{2a}$ is the sum of all weak
$(2a+1)$ block partitions
such that the zero block is always empty; so we have
upto $2a$ non-empty blocks.
This explains the term ``$2a$-shuffle''. 
If we assign labels $1,2,\ldots,a$ in decreasing order to the first 
$a$ blocks with the leftmost (signed) block labelled $a$, then
we get the random walk description for $S_{2a}$ that we started with. 

With the partition description it is immediate that 
$S_{2a} S_{2b} = S_{4ab}$. 
This gives us a multiplicative shuffle semigroup $\{S_{2a} : a > 0
\} \cup S_1$
contained in $k[\sigma_1] = k[S_2]$.
It now follows that 
$\sigma_0=1$,$\sigma_1$,\ldots,$\sigma_{n}$ is a basis for
$k[\sigma_1]$ and that it is a shuffle algebra.

\medskip
\noindent
{\bf The odd part}.

This case is completely analogous to the even case,
the only difference being that now we put no
restriction on the zero block.
We begin with the element $\sigma^{\prime}_1 =
\sigma_{s_1}+\cdots+\sigma_{s_{n}}$; that is,
$\sigma^{\prime}_1$ is the sum of all the vertices in $\Sigma$.  
We also
define $\sigma^{\prime}_j$ as the sum of all the simplices of rank $j$ in $\Sigma$
for $j = 1,\ldots,n$. 
In terms of partitions, $\sigma^{\prime}_{j}$ is the sum of all
$(2j+1)$ block partitions.
Note that $\sigma^{\prime}_n = \sigma_n$.
Also put
$\sigma^{\prime}_0 = \sigma_0 = 1$.
Next 
define shuffles $S_{2a+1}$ for $a \geq 0$ by
\begin{equation}\label{e:orshuffle}
S_{2a+1} = \sum_{j=0}^n \binom{a}{j} \sigma^{\prime}_{j}.
\end{equation}
To express $S_{2a+1}$ in words, we again use the notion of a
weak partition. 
Exactly as in the even case, it follows that $S_{2a+1}$ is the sum
of all weak $(2a+1)$ block partitions.
The random walk description is obtained by assigning 
labels $0,1,\ldots,a$ in decreasing order to the first 
$(a+1)$ blocks with the leftmost (signed) block labelled $a$.

With the partition description it is immediate that 
$S_{2a+1} S_{2b+1} = S_{(2a+1)(2b+1)}$. 
This gives us a multiplicative shuffle semigroup $\{S_{2a+1} : a \geq 0\}$
contained in $k[\sigma^{\prime}_1] = k[S_3]$.
It now follows that 
$\sigma^{\prime}_0=1$,$\sigma^{\prime}_1$,\ldots,$\sigma^{\prime}_{n}$ is a
basis for $k[\sigma^{\prime}_1]$ and that it is a shuffle algebra.

\medskip
\noindent
{\bf The shuffle double.}

We now put the even and odd parts together.
Define $S = \{S_{2a} : a \geq 1\} \cup \{S_{2a+1} : a \geq 0\}$.
It is immediate from the descriptions of $S_{2a}$ and $S_{2a+1}$ in
terms of partitions that 
$S_{a} S_{b} = S_{ab}$ for all $a,b \geq 1$ irrespective of parity. 
Thus we get a multiplicative shuffle semigroup $S$
as a spanning set in $A = k[\sigma_1,\sigma^{\prime}_1]$.
It follows that $A$ is a $2n$-dimensional algebra with basis  
$\sigma_0=1,\sigma_1,\sigma^{\prime}_1,\ldots,\sigma_{n-1},\sigma^{\prime}_{n-1},\sigma_{n}=\sigma^{\prime}_{n}$.
We call it a shuffle ``double'' algebra. 
It is not a shuffle algebra because the first three conditions are each
violated by a factor of two.
For example, in condition $(3)$, instead of a single generator $\sigma_1$,
we have two generators $\sigma_1$ and $\sigma^{\prime}_1$.

We now show that $A \iso k^{2n}$.
Following the example of the riffle shuffle of type $A_{n-1}$, we write
equation~\eqref{e:ershuffle} as 
\begin{equation}\label{e:1} 
S_{2a}=\sum_{i=1}^n (2a)^i e_i, \quad \text{where} \quad e_i = \sum_{j=i}^n
c_{ij}\frac{\sigma_{j}}{2^j j!}\,. 
\end{equation}
Here $c_{ij}$ is the coefficient of $x^i$ in the polynomial
$x(x-2)(x-4)\cdots(x-2(j-1))$.
It is $2^{j-i}$ times the $c_{ij}$ that occurred in the riffle shuffle
of type $A_{n-1}$. For notational convenience, we set $e_0=0$.
Similarly, we rewrite equation~\eqref{e:orshuffle} as
\begin{equation}\label{e:2} 
S_{2a+1}=\sum_{i=0}^n (2a+1)^i e^{\prime}_i, \quad \text{where} \quad e^{\prime}_i = \sum_{j=i}^n
c^{\prime}_{ij}\frac{\sigma^{\prime}_{j}}{2^j j!}\,. 
\end{equation}
Here $c^{\prime}_{ij}$ is the coefficient of $x^i$ in the polynomial
$(x-1)(x-3)\cdots(x-(2j-1))$.
Note that $e^{\prime}_n = e_n = \frac{\sigma_n}{2^n n!}$.
We next define two different kinds of characters of $S$. 
For any non-negative integer $i$, define
$\chi_i(S_{a}) = a^i$ for $a\geq1$.
Also define 
$\chi^{\prime}_i(S_{2a}) = 0$ for $a\geq1$ and $\chi^{\prime}_i(S_{2a+1}) = (2a+1)^i$
for $a\geq0$.
Now observe that for any $s \in S$ 
\begin{equation*} 
s = \sum_{i=1}^n \chi_i(s) e_i + \sum_{i=0}^n \chi^{\prime}_i(s) (e^{\prime}_i - e_i).
\end{equation*}
This can be checked separately for $s=S_{2a},S_{2a+1}$ using
equations~\eqref{e:1} and~\eqref{e:2}. 
The second summation actually goes only till $n-1$, since $e^{\prime}_n = e_n$.
Now apply Proposition~\ref{p:ss} to $S$ and $A$ along with 
the $2n$ characters $\chi_i$ for $i=1,\ldots,n$
and $\chi^{\prime}_i$ for $i=0,1,\ldots,n-1$
to conclude that $A \iso k^{2n}$. 

It is clear that $k[\sigma_1] \isoto k^{n+1}$, 
where the $n+1$ factors correspond to the idempotents 
$e^{\prime}_0,e_1,e_2,\ldots,e_n$. 
The isomorphism maps $S_{2a}$ to $(0,(2a)^1,(2a)^2,\ldots,(2a)^n)$.
To see that $k[\sigma^{\prime}_1] \iso k^{n+1}$, 
first note that $\chi^{\prime}_i(s) = \chi_i(s)$ for $s \in \{S_{2a+1} : a \geq
0\}$. 
Hence we need to lump together the factors corresponding to $e_i$ and
$e^{\prime}_i - e_i$. Thus we get $k[\sigma^{\prime}_1] \isoto k^{n+1}$, with
$e^{\prime}_0,e^{\prime}_1,e^{\prime}_2,\ldots,e^{\prime}_n$ as the primitive idempotents.
The isomorphism maps $S_{2a+1}$ to $((2a+1)^0,(2a+1)^1,(2a+1)^2,\ldots,(2a+1)^n)$.

\section{Examples of type $D_n$}

We consider two examples, the side shuffle and the riffle shuffle.
In this section, the term partition always refers to a partition of
type $D_n$; see Appendix~\ref{subs:d}. 
Elements of the label set $I$ will be written as $s_1, s_2, \ldots, s_{n-2},u,v$;
see Figure 2.3 in Appendix~\ref{subs:cd}.
We also recall that the random walk operates on 
an almost signed deck of cards or a deck of type $D_n$. 
It is a deck in
which every card, except the bottomost, is signed.
Also, one of the two basic algebras of Section~\ref{subs:drs} 
does not strictly satisfy
condition $(2)$ in the definition of a shuffle algebra. 
This is relevant to the remark in Section~\ref{subs:defn}.

\subsection{The side shuffle}
The element of interest is $\sigma = \sigma_{s_1}$.  In terms of
partitions, $\sigma_{s_1}$ is the sum of all three block partitions
such that the first (and hence the last) block is a singleton.  
The associated random walk on an almost signed deck of $n$ cards
consists of removing a card at random and replacing it on top with
either the same sign or its reverse both with equal probability.
If the bottomost unsigned card is chosen then we give it a sign and
then erase the sign on the new bottomost card.

In addition to $\sigma = \sigma_1$,
we define $\sigma_j = \sigma_{J_j}$ where $J_j = \{s_1,\ldots,s_j\}
\subseteq I$ for $j = 1,\ldots,n-2$. 
Note that the choice for $\sigma_{n-1}$ is not immediately obvious.
We set 
$\sigma_{n-1} = \sigma_{\{s_1,\ldots,s_{n-1},u,v\}}$.
The motivation becomes clearer from the following.
As partitions, $\sigma_{j}$ is the sum of all
$(2j+1)$ block partitions such that the first (and hence the last) $j$
blocks are singletons.
Observe that for $j=n-1$, this does give us the above definition of
$\sigma_{n-1}$.
We also set $\sigma_n = 2 \sigma_{n-1}$.
Define $a$-shuffles for $a \geq 0$ by 
\begin{equation*} 
S_a = \sum_{j=0}^n S(a,j) \sigma_j,
\end{equation*}
where $S(a,j)$ are the signed Stirling numbers defined in
Section~\ref{subs:atss}.
The associated random walks on a deck of type $D_n$ are as follows.

$\sigma_j :$ Choose $j$ cards at random. Move them to the top in a
random order. Now with equal probability either flip or do not flip
the sign of the chosen cards.  
If the bottomost unsigned card is one of the chosen cards then we give
it a sign and 
then erase the sign on the new bottomost card.

$S_a :$ Assign the sign $+$ to the bottomost card. Now do the side
shuffle $S_a$ of type $B_n$. Drop the sign of the bottomost card.

Observe that $\sigma_n = 2 \sigma_{n-1}$. We did not have any such
relation for the side shuffle of type $B_n$. It is interesting how
this relation emerges by just making the bottomost card unsigned.
Also observe that $S_a S_b = S_{a+b}$. At this point, it is clear that
the rest is identical to the two sided shuffle for type
$A_{n-1}$ and we omit it. 

\subsection{The riffle shuffle} \label{subs:drs}

The shuffles $S_a$ in this example can be described in the same way as
the inverse $a$-shuffles of type $B_{n}$ with a minor modification.  
This is because we are now operating on a deck of type $D_n$
rather than of type $B_n$.
The difference between the two is that in the former the bottomost
card is unsigned.
We now describe $S_a$, the inverse $a$-shuffle of type $D_n$.

$S_a :$ Assign the sign $+$ to the bottomost card. Now do an inverse
$a$-shuffle of type $B_n$. Drop the sign of the bottomost card.

It is useful to describe these shuffles in another equivalent way
using the analogy with the previous situation 
rather than the end result. For that,
we first define the inverse riffle shuffle of type $D_n$.

For every signed card, we either flip or do not flip its sign with equal
probability. The cards whose signs were flipped move below the
unsigned card (in the reverse relative order). The rest stay on top.

Note that after this shuffle, we may not have a deck of type $D_n$, 
since the unsigned card is not necessarily at the bottom. Hence we
define a \emph{correction} operation which assigns a sign to the unsigned
card and then erases the sign of the bottomost card.
We now give a description of our shuffles as operators on a deck of
type $D_n$.

$S_{2a}$: We do a usual inverse $a$-shuffle with labels
$1,2,\ldots,a$. Then within each block with a fixed label we do an
inverse riffle shuffle of type $D_n$ or $B_n$ depending on whether the
block contains the unsigned card or not. Lastly, we do the correction.

$S_{2a+1}$: We do a usual inverse $(a+1)$-shuffle with labels
$0,1,2,\ldots,a$. Then we repeat the above on each
block except the one labelled 0 and then do the correction.

Note that $S_2$ is just an inverse riffle shuffle of type $D_n$
followed by the correction.

\medskip
\noindent
{\bf The even part}.

The element of interest is $\sigma_1 = \sigma_{u} + \sigma_{v}$. 
In terms of partitions, $\sigma_{1}$ is the sum of all
three block partitions such that the central
block is  empty.

Put $\sigma_j = \sum_{\abs{J}=j,u \in J}\sigma_J + \sum_{\abs{J}=j,v \in J}\sigma_J$ for $j=1,2,\ldots,n$.
Explicitly for $n=3$, we have
$\sigma_1 = \sigma_{u} + \sigma_{v},
\sigma_2 = \sigma_{\{s,u\}} + \sigma_{\{s,v\}} + 2 \sigma_{\{u,v\}},
\sigma_3 = 2 \sigma_{\{s,u,v\}}$.
Note that we consider only those faces whose type contains either $u$
or $v$. If it contains both $u$ and $v$ then we put a coefficient of
$2$ in front.
In terms of partitions, $\sigma_{j}$ is the sum of all
$(2j+1)$ block partitions such that the central
block is  empty.
In this description we follow the convention that the partitions 
$(\{1,\overline{2}\},\{3\},\{\},\{\overline{3}\},\{2,\overline{1}\})$ 
and
$(\{1,\overline{2}\},\{\overline{3}\},\{\},\{3\},\{2,\overline{1}\})$
are counted separately, even though they are both equal to 
$(\{1,\overline{2}\},\{3,\overline{3}\},\{2,\overline{1}\})$.
This explains the coefficient of $2$ in front of faces whose type
contains both $u$ and $v$.

\medskip
\noindent
{\bf The odd part}.

The element of interest is $\sigma^{\prime}_1 = (\sigma_{u} + \sigma_{v}) +
(\sigma_{s_1} + \ldots + \sigma_{s_{n-2}}) + \sigma_{\{u,v\}}$. 
Note that in contrast to all our examples so far,
$\sigma^{\prime}_1$ is a combination of elements of $k\Sigma_1$
and $k\Sigma_2$. Next we define
$\sigma^{\prime}_j = \sigma_j + \sum_{\abs{J}=j,u \notin J, v \notin J}\sigma_J +
\sum_{\abs{J}=j+1,\{u,v\} \subset J}\sigma_J$ for $j=1,2,\ldots,n$.
Here $\sigma_j$ is as defined in the even case.
Thus we see that $\sigma^{\prime}_j$ is a combination of elements of $k\Sigma_j$
and $k\Sigma_{j+1}$.
This is the only example where the geometric description of 
objects of interest is somewhat complicated.
For better motivation we now consider partitions. 
We have already interpreted $\sigma_j$ in terms of partitions.
Note that the two summation terms that we added give the sum 
of all
$(2j+1)$ block partitions with a non-empty central
block.
Hence we see that $\sigma^{\prime}_j$ is the sum 
of all
$(2j+1)$ block partitions.
While counting partitions, we adopt the same convention as in the even case.

\medskip
\noindent
{\bf The shuffle double.}

We define shuffles $S_{2a}$ and $S_{2a+1}$ by
equations~\eqref{e:ershuffle} and~\eqref{e:orshuffle} and the analysis is
identical word for word to the case of $B_n$. We only point out that
partition now means partition of type $D_n$ and $\sigma_j$ and
$\sigma^{\prime}_j$ refer to the definitions that we made above. 
The shuffles can be described using weak partitions 
and this leads to the two random walk interpretations
that we gave earlier.
They are slightly complicated now
because of the more involved nature
of the partitions of type $D_n$.
A clear conceptual explanation of the close connection between the
riffle shuffles of type $B_n$ and $D_n$ is given in the next section.

\begin{remark}
The riffle shuffle examples in Sections~\ref{subs:ars}
and~\ref{subs:brs}(odd part) show that for types $A_{n-1}$ and $B_n$,
$\sigma_i \sigma_j = \sigma_j \sigma_i$, 
where
$\sigma_i$ is the sum of all the simplices of rank $i$.
This is because $\sigma_i$ and $\sigma_j$
are elements of a commutative (shuffle) algebra.
The fact that
$\sigma_i \sigma_j = \sigma_j \sigma_i$
can also be proven by a direct geometric argument. The key property is
that the simplicial 
complex induced on the support of any face of the Coxeter complex of type
$A_{n-1}$ or $B_n$ is again of the same type.
This property fails for $D_n$.
It is incomplete in this geometric sense.
Note that the $\sigma_i$'s as defined above did not play any role in the
riffle shuffle for $D_n$.
In fact, their role for type $D_n$ in this theory is far from being clear.
\end{remark}

\section{Maps between Coxeter complexes} \label{s:maps}

Recall that the reflection arrangements of type $B_n$, $D_n$
and $A_{n-1}$ are given by the hyperplanes
$x_i= \pm x_j, x_i = 0$; $x_i= \pm x_j$ and $x_i=x_j$ 
$(1\le i<j\le n)$ respectively.
Let $\Sigma(B_n)$, $\Sigma(D_n)$
and $\Sigma(A_{n-1})$ be the corresponding Coxeter complexes.
Observe that the arrangement for $D_n$ (resp. $A_{n-1}$) is
obtained from the one for $B_n$ (resp. $D_n$) by deleting some
hyperplanes. 
Following the discussion in Appendix~\ref{subs:map}, we get maps
$\Sigma(B_n) \rightarrow \Sigma(D_n) \rightarrow \Sigma(A_{n-1})$.
These are semigroup homomorphisms and
they induce algebra homomorphisms
$k\Sigma(B_n) \rightarrow k\Sigma(D_n) \rightarrow k\Sigma(A_{n-1})$.
We will use the language of
partitions to make these maps more explicit.

\begin{figure}[hbt]
\centering
\begin{tabular}{c@{\qquad}c}
\mbox{\epsfig{file=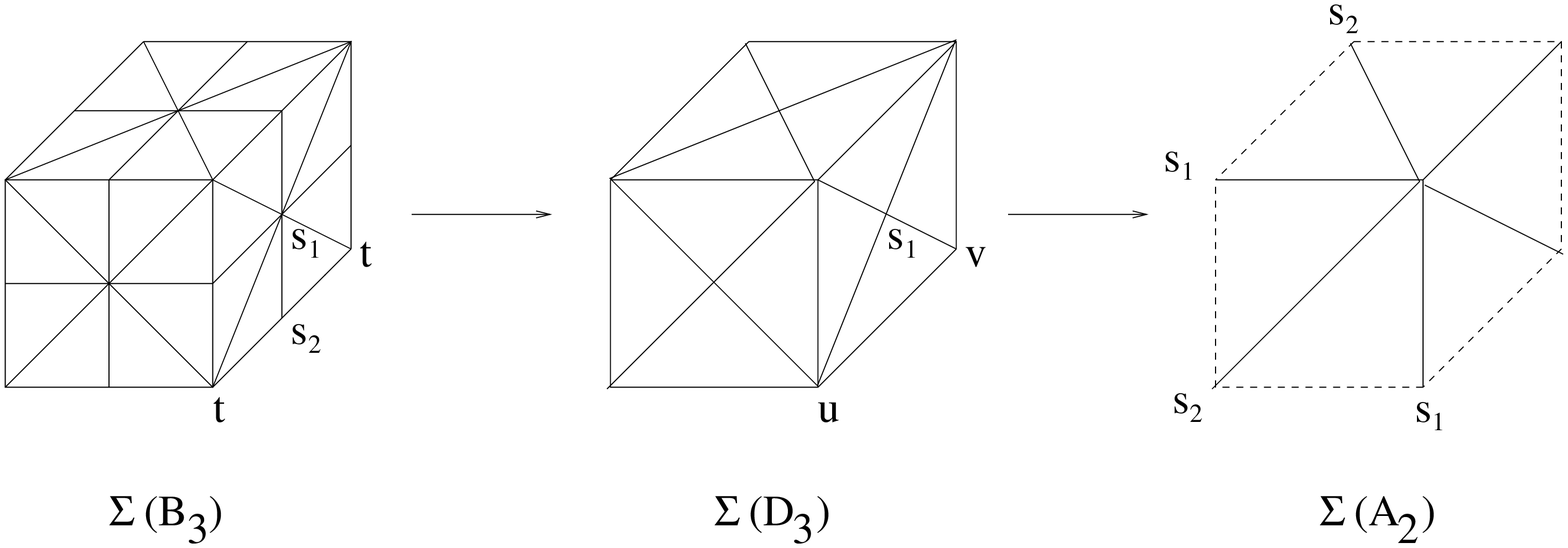,height=2.7cm,width = 12cm}}
\end{tabular}
\caption{The case $n=3$.}
\label{maps}
\end{figure}

\noindent
Figure~\ref{maps} shows the intersection of the hyperplane arrangements
for $B_3$, $D_3$ and $A_2$
with the boundary of a cube centred at the origin.
For $B_3$ and $D_3$, 
this gives us the corresponding Coxeter complex.
For $A_2$, we do not quite get $\Sigma(A_2)$ 
because the braid arrangement is not essential.
So we cut the braid arrangement 
by the hyperplane $x_1+x_2+x_3=0$.
The complex $\Sigma(A_2)$ can be seen in Figure~\ref{maps}
as the dotted hexagon 
with three vertices each of type $s_1$ and $s_2$.
Also in Figure~\ref{maps},
we have labelled only some of the vertices
since they uniquely determine the remaining labels.

\subsection{The map $\Sigma(B_n) \rightarrow \Sigma(D_n)$} \label{subs:map1}

For the case $B_n$ (resp. $D_n$) we have split the partition
description for a face $F$ into 4 (resp. 3) cases; see
Appendix~\ref{subs:b} and~\ref{subs:d}. From those descriptions, it is
straightforward to describe the map. 
Let $K \subset \{s_1,\ldots,s_{n-2}\}$. 

\begin{itemize}
\item[(i)]
A face of type $K$ in $\Sigma(B_n)$ maps to a face of type $K$ in
$\Sigma(D_n)$. 

\item[(ii)]
A face of type $K \cup t$ in $\Sigma(B_n)$  maps to a face either of type 
$K \cup u$ or $K \cup v$ depending on a parity condition. 

\item[(iii)]
A face of type $K \cup s_{n-1}$ in $\Sigma(B_n)$ maps to a face of one
higher dimension of type  
$K \cup \{u,v\}$. 
For $n=3$, this is also clear from Figure~\ref{maps},
where we see that 
$\{s_2\} \mapsto \{u,v\}$ and $\{s_1,s_2\}\mapsto\{s_1,u,v\}$. 

\item[(iv)]
A face of type $K \cup \{s_{n-1},t\}$ in $\Sigma(B_n)$ maps to a face of type 
$K \cup \{u,v\}$. 
\end{itemize}
The maps in the first three cases are bijective.
More precisely, the faces of type $K$ in $\Sigma(B_n)$
are in bijection with the faces of type $K$ in $\Sigma(D_n)$
and so on. However in the last case, 
the map is two to one. This is again
illustrated in Figure~\ref{maps},
where we see that there are two faces of type $\{s_1,s_2,t\}$
which map to a face of type $\{s_1,u,v\}$.

Recall that $\sigma_J$ is the sum of all faces of type $J$.
From the above discussion it follows that 
for the map $k\Sigma(B_n) \rightarrow k\Sigma(D_n)$,
we have
$\sigma_{K} \mapsto \sigma_{K}$,
$\sigma_{K \cup s_{n-1}} \mapsto \sigma_{K \cup \{u,v\}}$,
$\sigma_{K \cup t} \mapsto \sigma_{K \cup u} + \sigma_{K \cup v}$ and
$\sigma_{K \cup \{s_{n-1},t\}} \mapsto 2 \sigma_{K \cup \{u,v\}}$.

\subsection{The map $\Sigma(B_n) \rightarrow \Sigma(A_{n-1})$} \label{subs:map2}
It is easier to describe this composite map than the map from
$\Sigma(D_n)$ to $\Sigma(A_{n-1})$. In contrast to the previous map,
this composite map is hard to describe using face types but easy to
describe using partitions. 
To get an ordered partition of $[n]$, 
starting with an anti-symmetric
ordered partition of $[n,\overline{n}]$, we simply forget the set $[\overline{n}]$. 
For example,
$(\{2\},\{\overline{3}\},\{1,\overline{1}\},\{3\},\{\overline{2}\})$
maps to $(\{2\},\{1\},\{3\})$.
We explain some special cases.

Observe that a face of type $s_1$ maps to a face either of type $s_1$
or $s_{n-1}$ depending on whether the element in the singleton first
block has a plus sign or minus. It follows that
for the map $k\Sigma(B_n) \rightarrow k\Sigma(A_{n-1})$,
$\sigma_{s_1}$ maps to $\sigma_{s_1}+\sigma_{s_{n-1}}$.
This is the primary reason why the side shuffle of type $B_n$
is more closely related to the two-sided shuffle 
(rather than the side shuffle) of type $A_{n-1}$.

Next we show that inverse $a$-shuffles map to inverse $a$-shuffles.
For that
observe that weak $(2a+1)$ block partitions of type $B_n$
correspond to weak $(2a+1)$ block partitions of type $A_{n-1}$.
For example,
$(\{\overline{2}\},\{\},\{1,\overline{1}\},\{\},\{{2}\})
\leftrightarrow (\{\},\{\},\{1\},\{\},\{{2}\})$.
Similarly, weak $(2a+1)$ block partitions of type $B_n$ with an empty
zero block correspond to weak $2a$ block partitions of type $A_{n-1}$.
For example, when $a=1$, we get that 
$\sigma_{t}$ maps to $\sigma_{s_1}+\ldots+\sigma_{s_{n-1}}+2\sigma_0$.
In other words, 
the inverse riffle shuffle of type $B_n$ maps to the inverse
riffle shuffle of type $A_{n-1}$.

\subsection{Maps between side shuffles}

We recall the definitions of $\sigma_j$ for the side
shuffles of type $B_n$ and $D_n$ and the two-sided shuffle of type
$A_{n-1}$.

\begin{itemize}
\item[$B_n$] $: 
\sigma_j = \sigma_{J_j}$ where $J_j = \{s_1,\ldots,s_j\}
\subseteq I$ for $j = 1,\ldots,n$. 

\item[$D_n$] $: 
\sigma_j = \sigma_{J_j}$ where $J_j = \{s_1,\ldots,s_j\}
\subseteq I$ for $j = 1,\ldots,n-2$ and $J_{n-1} = \{s_1,\ldots,s_{n-1},u,v\}$. 
Also $\sigma_n = 2 \sigma_{n-1}$.

\item[$A_{n-1}$] $: 
\sigma_j = \sum_{k=0}^j \binom{j}{k}\sigma_{J_{j,k}}$ where
$J_{j,k} = \{s_1,\ldots,s_k\} \cup \{s_{n-(j-k)},\ldots,s_{n-1}\}
\subseteq I$ for $j = 1,\ldots,n-1$ and $\sigma_n = 2 \sigma_{n-1}$.

\end{itemize}
We claim that 
$\sigma_j \mapsto \sigma_j \mapsto \sigma_j$
under the maps 
$k\Sigma(B_n) \rightarrow k\Sigma(D_n) \rightarrow k\Sigma(A_{n-1})$.
From the discussion in Section~\ref{subs:map1},
it is easy to see the claim for the first map.
A point worth noting is that $\sigma_n \mapsto 2 \sigma_{n-1}$
and the claim holds for $j=n$
because of the relation $\sigma_n = 2 \sigma_{n-1}$ 
for $D_n$.
Similarly, using the discussion in Section~\ref{subs:map2},
we can show the claim for the composite map.
In earlier sections,
the introduction of $\sigma_n$ and the relation 
$\sigma_n = 2 \sigma_{n-1}$
for types $D_n$ and $A_{n-1}$ was somewhat artificial and justified 
only from the random walk perspective.
However we now also have a clear geometric perspective.

Recall that the shuffles are given by the same formula 
$S_a = \sum_{j=0}^n S(a,j) \sigma_j$
in all three cases with $S(a,j)$ being the signed Stirling numbers. 
Hence we get
$S_a \mapsto S_a \mapsto S_a$.
To summarise, we have:

\begin{displaymath}
\begin{array}{ccccc}
\left\{
\begin{array}{c}
\text {side shuffle}\\
\text {of type $B_n$}
\end{array}
\right\}&\longrightarrow&\left\{
\begin{array}{c}
\text {side shuffle}\\
\text {of type $D_n$}\\
\end{array}
\right\}&\isoto&\left\{
\begin{array}{c}
\text {two sided shuffle}\\
\text {of type $A_{n-1}$}\\
\end{array}
\right\}\\
\sigma_j & \mapsto & \sigma_j & \mapsto & \sigma_j \\
S_a & \mapsto & S_a & \mapsto & S_a 
\end{array}
\end{displaymath}

\subsection{Maps between riffle  shuffles}

Recall the definitions of $\sigma_j$ and $\sigma^{\prime}_j$ for the riffle
shuffles of type $B_n$ and $D_n$ from 
Sections~\ref{subs:brs} and \ref{subs:drs}. 
Following the discussion in
Section~\ref{subs:map1}, it is clear that 
$\Sigma(B_n) \rightarrow \Sigma(D_n)$
maps $\sigma_j$ to $\sigma_j$ and $\sigma^{\prime}_j$ to $\sigma^{\prime}_j$.
So the map restricted to the algebras $A = k[\sigma_1,\sigma^{\prime}_1]$ 
is an algebra map that maps a basis to a basis.
Hence it is an isomorphism. 
It also clearly maps $S_a$ to $S_a$.
Again we see that the somewhat unmotivated definitions of
$\sigma_j$ and $\sigma^{\prime}_j$ for $D_n$
have a more natural meaning as the images of the corresponding 
elements in $\Sigma(B_n)$.

For the map $\Sigma(B_n) \rightarrow \Sigma(A_{n-1})$, 
the discussion at the end of Section~\ref{subs:map2} 
and the interpretation of the inverse $a$-shuffles $S_a$ as weak partitions shows
that $S_a \mapsto S_a$. This gives us a surjective map
$A = k[\sigma_1,\sigma^{\prime}_1] = k[S_2,S_3] \rightarrow A = k[S_2]$.  
However we do not know a good way to describe 
the images of $\sigma_j$ and $\sigma^{\prime}_j$ under this map. 
To summarise, we have:

\begin{displaymath}
\begin{array}{ccccc}
\left\{
\begin{array}{c}
\text {riffle shuffle}\\
\text {of type $B_n$}
\end{array}
\right\}&\isoto&\left\{
\begin{array}{c}
\text {riffle shuffle}\\
\text {of type $D_n$}\\
\end{array}
\right\}&\longrightarrow&\left\{
\begin{array}{c}
\text {riffle shuffle}\\
\text {of type $A_{n-1}$}\\
\end{array}
\right\}\\
\sigma_j & \mapsto & \sigma_j & \mapsto & ? \\
\sigma^{\prime}_j & \mapsto & \sigma^{\prime}_j & \mapsto & ? \\
S_a & \mapsto & S_a & \mapsto & S_a 
\end{array}
\end{displaymath}

\vanish{
\matrix
{\left\{
\begin{tabular}{c}
side shuffle\\
of type $B_n$\\
\end{tabular}
\right\}&\rightarrow&\left\{
\begin{tabular}
side shuffle\\
of type $B_n$\\
\end{tabular}
\right\}&\rightarrow&\left\{
\begin{tabular}
side shuffle\\
of type $B_n$\\
\end{tabular}
\right\}\cr}
}

\begin{remark}
The Coxeter complex of type $B_n$ is probably the correct place
to look for shuffle algebras.
The images of these shuffle algebras under our maps
would then yield shuffle algebras (in a suitably generalised sense)
of type $D_n$ and $A_{n-1}$.
Among the examples in this paper,
the only one that we failed to describe in this manner
was the side shuffle of type $A_{n-1}$.
\end{remark}

\section{Generalisation to buildings}\label{s:q}

We first give a brief review of buildings.
For the general theory
see~\cite{\brown,\schar,\tits}. 
Let $W$ be a Coxeter group and $\Sigma(W)$ its Coxeter complex.
Roughly,
a building $\Delta$ of type $W$ is a union of subcomplexes $\Sigma$ (called
apartments) which fit together nicely.
Each apartment $\Sigma$ is isomorphic to $\Sigma(W)$.
As a simplicial complex, $\Delta$ is pure and labelled.
The term \emph{pure} means that all maximal simplices (chambers)
have the same dimension. 
For any two simplices in $\Delta$, there is an apartment $\Sigma$
containing both of them.
Using this fact, we can define a product on $\Delta$ as follows.

For $x,y \in \Delta$, we choose an apartment $\Sigma$ containing $x$ and $y$
and define $xy$ to be their product in $\Sigma$.
Since $\Sigma$ is a Coxeter complex, we know how to do this.
Furthermore, it can be shown that the product does not depend on the
choice of $\Sigma$.
So this defines a product on $\Delta$ 
with the set of chambers $\C$ contained as an ideal.
As in Section~\ref{subs:rw},
we can now define a random walk on $\C$ 
starting with a probability distribution on $\Delta$.
Unfortunately, the product on $\Delta$ is not associative,
that is,
$\Delta$ is no longer a semigroup.
As a result $k\C$ is not a module over $k\Delta$
and so our algebraic methods break down.
However, we know three examples where they do work.
This is because in these cases, the subalgebra of interest,
namely $k[\sigma]$, turns out to be associative.

A theorem of Tits \cite[Theorem 11.4]{\tits} says that
if $\Delta$ is a finite, irreducible Moufang building
then it is the building of an
absolutely simple algebraic group $G$ over a finite field $\fq$.
The buildings we consider satisfy this hypothesis;
hence in each example we will also specify the algebraic group $G$.
The group $G$ acts by simplicial type-preserving automorphisms on
$\Delta$. To avoid getting into details, 
we just mention that the action of $G$ is very closely related to
the geometry of the building.

The three examples we consider can be thought of as 
$q$-analogues of the side shuffles of 
type $A_{n-1}$, $B_n$ and $D_n$.
We use the same definitions for the $\sigma_j$'s as before,
except that we now apply them to the building $\Delta$
rather than the Coxeter complex $\Sigma$.
We will use the notation $[j] = 1 + q + \ldots + q^{j-1}$
to denote the $q$-numbers. They will show up a lot in our analysis.

\subsection{The q-side shuffle for $A_{n-1}$} \label{subs:qass}

We first briefly describe the building of type $A_{n-1}$
associated to the algebraic group $GL_n(\fq)$ (also $SL_n(\fq)$).
Let $V$ be the $n$-dimensional vector space over $\fq$
and $\mathcal{L}_{n}$ be the lattice of subspaces
of $V$. 
The building of type $A_{n-1}$ associated to $GL_n(\fq)$ 
is simply the flag (order) complex
$\Delta(\mathcal{L}_n)$. 
It is a labelled simplicial complex.
A vertex of $\Delta(\mathcal{L}_n)$ is by definition a proper subspace of $V$
and we label it by its dimension.
To be consistent with earlier notation,
we label the vertices $s_1, \ldots, s_{n-1}$ 
instead of just $1, \ldots, n-1$.
In particular, vertices of type $s_1$ are the one dimensional
subspaces of $V$.

Let $\mathcal{B}_{n}$ be
the Boolean lattice of rank $n$ consisting
of all subsets of an $n$-set ordered under inclusion.
Note that the order complex $\Delta(\mathcal{B}_n)$ is the Coxeter complex of type
$A_{n-1}$, 
which we earlier denoted by $\Sigma(A_{n-1})$.
Also note that a
choice of a basis for $V$ gives an 
embedding of $\mathcal{B}_n$ into $\mathcal{L}_n$.
The subcomplexes $\Delta(\mathcal{B}_n)$, for various embeddings $\mathcal{B}_n \into \mathcal{L}_n$,
play the role of apartments.
Alternatively, it is also useful to think of $\mathcal{B}_n$ as a specialisation of $\mathcal{L}_n$
for the degenerate case $q=1$.
In this sense $\Delta(\mathcal{L}_n)$
may be regarded as a $q$-analogue of $\Delta(\mathcal{B}_n) = \Sigma(A_{n-1})$.
For the rest of the section
we will denote the building $\Delta(\mathcal{L}_n)$ by $\Delta$.

The product in $\Delta$ can be made more explicit. 
A face of $\Delta$ is a chain in $\mathcal{L}_n$
and a chamber is a maximal chain in $\mathcal{L}_n$.
For $x, y \in \Delta$, the face $xy$ is the chain in $\mathcal{L}_n$ obtained by
refining the chain $x$ by the chain $y$ (using meets and joins as in a
Jordan--H\"{o}lder product).
This is a generalisation of the product in $\Sigma(A_{n-1})$
which we defined
in terms of refinement of one ordered partition by another.

The side shuffle of type $A_{n-1}$ (Section~\ref{subs:ass})
generalises in a straightforward way to this setting. Let
$\sigma_1$ be the sum of all the vertices of type $s_1$,
$\sigma_2$ be the sum of all the edges of type $s_1 s_2$ and so on till
$\sigma_{n-1}$ which is the sum of all chambers of $\Delta$.
The number of summands in $\sigma_{j}$, 
which is same as the number of faces of type
$s_1 s_2 \ldots s_j$, is
$[n][n-1] \ldots [n-j+1]$.
It follows directly from the definition that 
$\sigma_j \sigma_1 = \sigma_1 \sigma_j = 
[j] \sigma_j + q^j \sigma_{j+1}$ if $j < n-1$ and 
$\sigma_{n-1} \sigma_1 = [n] \sigma_{n-1}$.
As a consistency check, 
we can also verify that the number of terms 
on each side is the same.
These facts imply that $\sigma_1$ is power associative.
Hence we may write the following without ambiguity. 
\begin{equation}\label{e:qcob}
q^{\binom{j+1}{2}}\sigma_{j+1} = \sigma_1(\sigma_1-[1])\ldots(\sigma_1-[j]).
\end{equation}
This along with
$\sigma_{n-1} \sigma_1 = [n] \sigma_{n-1}$ implies that both
$\sigma_0=1$,$\sigma_1$,\ldots,$\sigma_{n-1}$ and
$\sigma_1^0=1$,$\sigma_1^1$,\ldots,$\sigma_1^{n-1}$
form a basis for the associative algebra $A = k[\sigma_1]$.
The basic relation satisfied by $\sigma_1$ is 
$\sigma_1(\sigma_1-[1])\ldots(\sigma_1-[n-2])(\sigma_1-[n]) = 0$.
Hence
$A \isoto \frac{k[x]}{x(x-[1])\ldots(x-[n-2])(x-[n])}$,
where $\sigma_1$ maps to $x$.
Since the roots are distinct
we conclude that $A$ is split semisimple.

As before,
we extend the definition of $\sigma_j$ to any $j$
using~\eqref{e:qcob} to obtain 
$\sigma_n = \sigma_{n-1}$ and $\sigma_j = 0$ for
$j > n$.
Inverting equation~\eqref{e:qcob} formally, we get
\begin{equation*}
\sigma_1^a = \sum_{j=0}^a S(a,j) \sigma_j
\end{equation*}
where $S(a,j)$ satisfies the recursion:
$S(a,j) = [j] S(a-1,j) + q^{j-1} S(a-1,j-1)$ with $S(a,1) = 1$ and 
$S(a,a) = q^{\binom{a}{2}}$. 
The numbers $S(a,j)$ give a $q$-analogue to the Stirling numbers of
the second kind. For more information on these numbers, see~\cite{\wachswhite} and the references therein.

\begin{remark}
The analysis we gave for this example was the analogue of the
direct method of Section~\ref{subs:ass}.
The difficulty with the shuffle method is that
we do not know of a good random walk interpretation of $\sigma_1^a$.
However as the direct method shows,
the algebra $A$ satisfies all the properties of 
an additive shuffle algebra 
(except that now we are in the more general context of buildings).
There are also natural ways to define the q-analogue of the two-sided
shuffle and the riffle shuffle. However the analysis breaks down
because the basic elements of interest are not
power associative.
\end{remark}

\subsection{The q-side shuffle for $B_{n}$} \label{subs:qbss}

There are two possibilities for the building depending 
on whether the associated algebraic group is
the symplectic group $Sp_{2n}(\fq)$ or
the orthogonal group $O_{2n}(\fq)$.
We consider them separately.
They may both be regarded as $q$-analogues of $\Sigma(B_n)$.
The side shuffle for $B_{n}$ generalises to both cases.
As expected, the computations have the same spirit as the previous section. 

\medskip
\noindent
{\bf The symplectic case.}
Let $V$ be a (even dimensional) vector space over $\fq$ 
with a skew-symmetric non-degenerate bilinear form $Q$.
More explicitly, $V$ has a basis
$e_1, e_2, \ldots, e_n, f_1, f_2, \ldots, f_n$ 
such that $Q(e_i,f_i) = 1$, $Q(f_i,e_i) = -1$
and all other pairings between basis vectors are zero.
The symplectic group $Sp_{2n}(\fq)$ 
is the group of automorphisms of $V$ 
that preserve the bilinear form $Q$.

The building of type $B_{n}$
associated to the symplectic group
is the flag (order) complex of all isotropic subspaces of $V$.
We will denote it by $\Delta$.
The definition of the product in $\Delta$ is similar to that in
Section~\ref{subs:qass}.
The difference occurs due to the fact that 
the sum of two isotropic subspaces may not be isotropic.

Let
$\sigma_1$ be the sum of all the vertices of type $s_1$,
$\sigma_2$ be the sum of all the edges of type $s_1 s_2$ and so on till
$\sigma_{n}$ which is the sum of all chambers of $\Delta$.
It follows from the definition that 
$\sigma_j \sigma_1 = \sigma_1 \sigma_j = 
(1+q^{2n-j})[j] 
\sigma_j + q^j \sigma_{j+1}$ if $j \leq n-1$ and 
$\sigma_{n} \sigma_1 = \sigma_{1} \sigma_n = [2n] \sigma_{n}$.
An easy way to check this is to first figure out the coefficient of 
$\sigma_{j+1}$. It is $q^j$, the same answer that we got for the
$A_{n-1}$ case.
There is no difference in the two situations
with respect to this coefficient.
With this information, we can find the coefficient of
$\sigma_{j}$ by simply counting the number of terms involved 
on both sides. 
The number of summands in $\sigma_{j}$, 
which is same as the number of faces of type
$s_1 s_2 \ldots s_j$, is
$[2n][2n-2] \ldots [2n-2j+2]$.
The identity is then a consequence of the simple relation
$[2n] = (1+q^{2n-j})[j] + q^j [2n-2j]$.
It is a good exercise to compute the coefficient of $\sigma_j$ directly.

These facts imply that $\sigma_1$ is power associative.
Hence we may write the following without ambiguity. 
\begin{equation}\label{e:qcob2}
q^{\binom{j+1}{2}}\sigma_{j+1} = \sigma_1(\sigma_1-(1+q^{2n-1})[1])\ldots(\sigma_1-(1+q^{2n-j})[j]).
\end{equation}
This along with
$\sigma_{n} \sigma_1 = [2n] \sigma_{n}$ implies that both
$\sigma_0=1$,$\sigma_1$,\ldots,$\sigma_{n}$ and
$\sigma_1^0=1$,$\sigma_1^1$,\ldots,$\sigma_1^{n}$
form a basis for the associative algebra $A = k[\sigma_1]$.
The basic relation satisfied by $\sigma_1$ is 
$\sigma_1(\sigma_1- (1+q^{2n-1})[1])\ldots(\sigma_1-(1+q^{n+1})[n-1])(\sigma_1-[2n]) = 0$.
Hence
$A \isoto \frac{k[x]}{x(x-(1+q^{2n-1})[1])\ldots(x-(1+q^{n+1})[n-1])(x-[2n])}$.
This shows that $A$ is split semisimple.

As before,
we extend the definition of $\sigma_j$ to any $j$
using~\eqref{e:qcob2} to obtain 
$\sigma_j = 0$ for
$j > n$.
Inverting equation~\eqref{e:qcob2} formally, we get
\begin{equation*}
\sigma_1^a = \sum_{j=0}^a S(a,j) \sigma_j
\end{equation*}
where $S(a,j)$ satisfies the recursion:
$S(a,j) = (1+q^{2n-j})[j] S(a-1,j) + q^{j-1} S(a-1,j-1)$ 
with $S(a,1) = (1+q^{2n-1})^{a-1}$ and 
$S(a,a) = q^{\binom{a}{2}}$. 
The numbers $S(a,j)$ give a $q$-analogue to the signed
Stirling numbers.
Note that they now depend on $n$.
We do not know whether they have been considered before or
whether they have an explicit formula.

\medskip
\noindent
{\bf The orthogonal case.}
Let $V$ be a (even dimensional) vector space over $\fq$ 
with a symmetric non-degenerate bilinear form $Q$.
We also assume that $V$ has an isotropic subspace of dimension $n$.
More explicitly, $V$ has a basis
$e_1, e_2, \ldots, e_n, f_1, f_2, \ldots, f_n$ 
such that $Q(e_i,f_i) = Q(f_i,e_i) = 1$
and all other pairings between basis vectors are zero.
The orthogonal group $O_{2n}(\fq)$ 
is the group of automorphisms of $V$ 
that preserve the bilinear form $Q$.

The building of type $B_{n}$
associated to the orthogonal group
is the flag (order) complex of all isotropic subspaces of $V$.
We will denote it by $\Delta(B_n)$.
The analysis parallels the symplectic case.
With the same definitions of the $\sigma_j$'s,
we get
$\sigma_j \sigma_1 = \sigma_1 \sigma_j = 
(1+q^{2n-j-1})[j] 
\sigma_j + q^j \sigma_{j+1}$ if $j \leq n-1$ and 
$\sigma_{n} \sigma_1 = \sigma_{1} \sigma_{n} = (1+q^{n-1})[n] \sigma_{n}$.
The reason for the difference is that now
the number of summands in $\sigma_{j}$ is
$(1+q^{n-1})[n](1+q^{n-2})[n-1] \ldots (1+q^{n-j})[n-j+1]$.
The rest is similar to the symplectic case.
We only point out that the Stirling numbers now satisfy the recursion:
$S(a,j) = (1+q^{2n-j-1})[j] S(a-1,j) + q^{j-1} S(a-1,j-1)$ 
with $S(a,1) = (1+q^{2n-2})^{a-1}$ and 
$S(a,a) = q^{\binom{a}{2}}$. 

\begin{remark}
The building $\Delta(B_n)$ that we considered here is not thick.
Classically, the group $O_{2n}(\fq)$ is considered 
to be of type $D_n$ since there is a thick building $\Delta(D_n)$
associated to it. We will discuss this next.
For more information on this, see~\cite[pgs 123-127]{\brown}.
\end{remark}

\subsection{The $q$-side shuffle for $D_n$} \label{subs:qdss}

Let $V$ and $Q$ be as in the orthogonal case of Section~\ref{subs:qbss}.
The building of type $D_n$,
which we denote by $\Delta(D_n)$,
is the flag complex of the following so-called ``oriflamme geometry''.
The vertices in $\Delta(D_n)$ are the non-zero isotropic subspaces of $V$
of dimension $\not=n-1$.
Two such subspaces are called incident 
if one is contained in the other or
if both have dimension $n$ and their intersection has dimension $n-1$.
A simplex in $\Delta(D_n)$ is a set of pairwise incident vertices.

There is a natural product preserving map 
$\Delta(B_n) \rightarrow \Delta(D_n)$,
which generalises the map $\Sigma(B_n) \rightarrow \Sigma(D_n)$
in Section~\ref{subs:map1}.
To describe this map,
we identify the vector space $V$ and the bilinear form $Q$
in the two cases.
A vertex in $\Delta(B_n)$ is a non-zero isotropic subspace of $V$.
If it has dimension $\not = n-1$
then we map it to itself.
Otherwise there are exactly two maximal isotropic subspaces 
of dimension $n$ in $V$ that contain it.
And we map it to the edge joining the two vertices in $\Delta(D_n)$
representing these two subspaces.
This again explains how a vertex of type $s_{n-1}$ 
maps to an edge of type $uv$. 

In the orthogonal case of Section~\ref{subs:qbss},
we have the shuffle algebra $k[\sigma_1]$ in $\Delta(B_n)$.
The image of this shuffle algebra under the above map
yields a shuffle algebra in $\Delta(D_n)$.
As a $q$-analogue of the map between side shuffles of 
type $B_n$ and $D_n$,
we get a map
\begin{displaymath}
\begin{array}{ccc}
\left\{
\begin{array}{c}
\text {$q$-side shuffle}\\
\text {of type $B_n$}
\end{array}
\right\}&\longrightarrow&\left\{
\begin{array}{c}
\text {$q$-side shuffle}\\
\text {of type $D_n$}\\
\end{array}
\right\}\\
\sigma_j & \mapsto & \sigma_j \\
\end{array}
\end{displaymath}
The $\sigma_j$'s on the right are defined 
as for the side shuffle of type $D_n$ and 
we again get the relation $\sigma_n = 2 \sigma_{n-1}$.
The basic relation satisfied by $\sigma_1$ is 
$\sigma_1(\sigma_1- (1+q^{2n-2})[1])\ldots(\sigma_1-(1+q^{n+1})[n-2])(\sigma_1-(1+q^{n-1})[n]) = 0$.
Note that the term corresponding to $n-1$ is absent.
The rest of the analysis is routine and we omit it.

\begin{remark}
In light of the discussion in Section~\ref{s:maps},
we may say the following.
The failure to get a $q$-two sided shuffle of type $A_{n-1}$
is the result of our failure to define a product preserving map
$\Delta(D_n) \rightarrow \Delta(A_{n-1})$.
Also we failed to get a $q$-riffle shuffle of type $B_n$
and as a result also failed on the other two fronts.
\end{remark}

\vanish{
The element of interest is $\sigma = \sigma_{s_1}$. 
In addition to $\sigma = \sigma_1$,
we define $\sigma_j = \sigma_{J_j}$ where $J_j = \{s_1,\ldots,s_j\}
\subseteq I$ for $j = 1,\ldots,n-2$. 
Next we set 
$\sigma_{n-1} = \sigma_{\{s_1,\ldots,s_{n-1},u,v\}}$, that is,
$\sigma_{n-1}$ is the sum of all chambers of $\Delta$.
It turns out that
Hence we may write the following without ambiguity. 
\begin{equation}\label{e:qcob3}
q^{\binom{j+1}{2}}\sigma_{j+1} = \sigma_1(\sigma_1-(1+q^{2n-2})[1])\ldots(\sigma_1-(1+q^{2n-j-1})[j]).
\end{equation}
This along with
$\sigma_{n-1} \sigma_1 = (1+q^{n-1})[n] \sigma_{n-1}$,
implies that both
$\sigma_0=1$,$\sigma_1$,\ldots,$\sigma_{n-1}$ and
$\sigma_1^0=1$,$\sigma_1^1$,\ldots,$\sigma_1^{n-1}$
form a basis for the associative algebra $A = k[\sigma_1]$.
We conclude that
$A \isoto \frac{k[x]}{x(x-(1+q^{2n-2})[1])\ldots(x-(1+q^{n+1})[n-2])(x-(1+q^{n-1})[n])}$.
This shows that $A$ is split semisimple.

As before,
we extend the definition of $\sigma_j$ to any $j$
using~\eqref{e:qcob} to obtain 

$\sigma_j = 0$ for
$j > n$.
Inverting equation~\eqref{e:qcob3} formally, we get
\begin{equation*}
\sigma_1^a = \sum_{j=0}^a S(a,j) \sigma_j
\end{equation*}
where $S(a,j)$ satisfies 
The numbers $S(a,j)$ give another $q$-analogue to the signed
Stirling numbers different from the one in the previous section.
}

\section{Future prospects}

We conclude by suggesting some problems for future consideration.

\subsection*{Shuffle algebras}
One problem is to classify them 
with minor modifications of the definition if necessary
(see remark in Section~\ref{subs:defn}). 
Also it is natural to ask whether there is a sensible
generalisation to the case of infinite Coxeter groups.

\subsection*{Buildings}
There is no good theory for random walks on buildings 
the major difficulty being the non-associativity of the product. 
Our methods in Section~\ref{s:q} suggest that it would be worthwhile
to systematically understand 
the power associative elements of a building
and also to classify product preserving maps between buildings.

\subsection*{Complex reflection groups}
In this paper, we considered only real reflection groups.
The difficulty in passing to the complex case 
is the absence of an analogue of the Coxeter complex.
As an example, consider the complex reflection group
$S_n \ltimes \Z_r^n$.
For $r=2$, it specialises to the Coxeter group of type $B_n$. 
The first step would be to generalise the semigroup
of ordered partitions of type $B_n$. 
However, it is not clear how to do this.

\subsection*{Multiplicities and derangement numbers}
In \cite{\ken}, Brown defines a random walk associated to matroids.
And he relates the eigenvalue multiplicities to
invariants of the lattice of flats.
He calls these invariants the generalised derangement numbers.
The motivating examples are the side and $q$-side shuffle of type $A_{n-1}$
and the multiplicities then are the usual and $q$-derangement numbers.
In this sense, we may think of ordinary matroids as related to
the Coxeter group of type $A_{n-1}$.

There is a notion of a $W$-matroid~\cite{\borovik} for any 
Coxeter group $W$. We ask whether it is possible to generalise the above 
to this setting.
As a positive result in this direction,
Bidigare \cite[pgs 147-148]{\bid} shows that 
for the side shuffle of type $B_n$,
the multiplicities are the signed derangement numbers.

\vanish{

Expect Tsetlin, riffle shuffle, etc.
}



\appendix

\section{The hyperplane face semigroup} \label{app:hyperplane}

Most of this section and a part of the next is taken directly from~\cite{\ken}.
More details concerning the material reviewed here can be found in
\cite{\bhr,\bbd,\red,\brown,\bd,\ot,\zie}.  Throughout this section
$\A=\{H_i\}_{i\in I}$ denotes a finite set of affine hyperplanes in
$V=\R^n$.  Let $H_i^+$ and $H_i^-$ be the two open halfspaces
determined by $H_i$; the choice of which one to call $H_i^+$ is
arbitrary but fixed.

\subsection{Faces and chambers}

The hyperplanes $H_i$ induce a partition of $V$ into convex sets
called \emph{faces} (or \emph{relatively open faces}).  These are the
nonempty sets $F \subseteq V$ of the form
\[
F = \bigcap_{i\in I} H_i^{\varepsilon_i},
\]
where $\varepsilon_i \in \{+,-,0\}$ and $H^0_i = H_i$.  Equivalently, if we
choose for each $i$ an affine function $f_i\colon V \rightarrow \R$
such that $H_i$ is defined by $f_i = 0$, then a face is a nonempty set
defined by equalities and inequalities of the form $f_i > 0$, $f_i <
0$, or $f_i = 0$, one for each $i \in I$.  The sequence $\varepsilon =
(\varepsilon_i)_{i\in I}$ that encodes the definition of $F$ is called the
\emph{sign sequence} of $F$ and is denoted $\varepsilon(F)$.

The faces such that $\varepsilon_i \ne 0$ for all $i$ are called
\emph{chambers}.  They are convex open sets that partition the
complement $V - \bigcup_{i\in I} H_i$.  In general, a face $F$ is
open relative to its \emph{support}, which is defined to be the affine
subspace
\[
\supp F = \bigcap_{\varepsilon_i (F) = 0} H_i.
\]
Since $F$ is open in $\supp F$, we can also describe $\supp F$ as the
affine span of~$F$.

\subsection{The face relation} 

The \emph{face poset} of $\A$ is the set $\F$ of faces, ordered as
follows:  $F \le G$ if for each $i \in I$ either $\varepsilon_i(F) = 0$ or
$\varepsilon_i(F) = \varepsilon_i(G)$.  In other words, the description of $F$
by linear equalities and inequalities is obtained from that of $G$ by
changing zero or more inequalities to equalities.
We say that $F$ is a face of $G$.
Note that the chambers are precisely the maximal elements 
of the face poset.
We denote the set of chambers by $\C$.

\subsection{Product}
The set $\F$ of faces is also a semigroup.  Given $F,G \in \F$, their
\emph{product} $FG$ is the face with sign sequence
\[
\varepsilon_i(FG) = 
\begin{cases}
  \varepsilon_i(F)  &\text{if $\varepsilon_i(F) \ne 0$} \\
         \varepsilon_i(G) &\text{if $\varepsilon_i(F) = 0$.}
\end{cases}
\]
This has a geometric interpretation:  If we move on a straight line
from a point of~$F$ toward a point of $G$, then $FG$ is the face we
are in after moving a small positive distance.  
Hence, we may think of $FG$ as the projection of $G$ on $F$.
Notice that the face
relation can be described in terms of the product:  One has
\begin{equation} \label{e:poset}
F\le G \iff FG=G.
\end{equation}
Also note that the set of chambers $\C$ is an ideal of $\F$.

\subsection{The semilattice of flats}

A second poset associated with the arrangement $\A$ is the
\emph{semilattice of flats}, also called the \emph{intersection
semilattice}, which we denote by~$\L$.  It consists of all nonempty
affine subspaces $X \subseteq V$ of the form $X = \bigcap_{H\in\A'}
H$, where $\A' \subseteq \A$ is an arbitrary subset (possibly empty).
We order $\L$ by inclusion.  [Warning: Many authors order $\L$ by
reverse inclusion.]  Notice that any two elements $X,Y$ have a
least upper bound $X\vee Y$ in~$\L$, which is the intersection of
all hyperplanes $H\in\A$ containing both $X$ and~$Y$; hence $\L$
is an \emph{upper semilattice} (poset with least upper bounds).  It is
a lattice if the arrangement $\A$ is \emph{central}, i.e., if
$\bigcap_{H\in\A} H \ne\emptyset$.  Indeed, this intersection is then
the smallest element of~$\L$, and a finite upper semilattice with a
smallest element is a lattice \cite[Section 3.3]{\stanley}.  The
support map gives a surjection
\[
\supp\colon\F\onto\L,
\]
which preserves order and also behaves nicely with respect to the
semigroup structure.  Namely, we have
\begin{equation}
\supp(FG) = \supp F \vee \supp G
\end{equation}
and
\begin{equation}
FG = F \iff \supp G \le \supp F.
\end{equation}

\subsection{The forgetful map} \label{subs:map}
Let $\A=\{H_i\}_{i\in I}$ and $\A^{\prime}=\{H_i\}_{i\in I^{\prime}}$ be two hyperplane
arrangements such that $I^{\prime} \subset I$, that is, the arrangement $\A^{\prime}$ is
obtained from $\A$ by deleting some of the hyperplanes.
Let $\F(A)$ and $\F(A^{\prime})$ be the respective hyperplane face
semigroups. Then there is a natural map $\F(A) \rightarrow \F(A^{\prime})$.
An element of $\F(A)$ is a face $F$ with a sign sequence, say 
$\varepsilon = (\varepsilon_i)_{i\in I}$. We map it to the face of
$\A^{\prime}$ with the sign sequence $\varepsilon^{\prime} = (\varepsilon_i)_{i\in I^{\prime}}$.
In other words, we forget the signs of the elements of $I$ that are not in
$I^{\prime}$.
It is clear that this map is a semigroup homomorphism and preserves
the face relation. The second fact is implied by the first in view of
equation~\eqref{e:poset}. 

\subsection{Spherical representation} \label{sub:cell}

Suppose now that $\A$ is a \emph{central} arrangement, i.e., that the
hyperplanes have a nonempty intersection.  We may assume that this
intersection contains the origin.  Suppose further that $\bigcap_{i\in
I} H_i = \{0\}$, in which case $\A$ is said to be \emph{essential}.
(There is no loss of generality in making this assumption; for if it
fails, then we can replace $V$ by the quotient space $V/\bigcap_i H_i$.)
The hyperplanes then induce a cell-decomposition of the unit sphere,
the cells being the intersections with the sphere of the faces
$F\in\F$.  Thus $\F$, as a poset, can be identified with the poset of
cells of a regular cell-complex $\Sigma$, homeomorphic to a sphere.
Note that the face $F=\{0\}$, which is the identity of the
semigroup~$\F$, is not visible in the spherical picture; it
corresponds to the empty cell.  The cell-complex $\Sigma$ plays a
crucial role in~\cite{\bd}, to which we refer for more details.

Thus the cell-complex $\Sigma$ is a semigroup 
which we call the hyperplane face semigroup.
The maximal cells (which we again call chambers and denote $\C$)
is an ideal in $\Sigma$.

\section{Reflection arrangements} \label{app:Coxeter}

We work with an
arbitrary finite Coxeter group~$W$ and its associated hyperplane face
semigroup $\Sigma$ (the Coxeter complex of~$W$).  
But we will 
explain everything in concrete terms for the cases $W=A_{n-1},B_n,D_n$.
This should make the discussion accessible to readers unfamiliar with
Coxeter groups.

\subsection{Finite reflection groups} \label{sub:reflection}

We begin with a very quick review of the basic facts that we need
about finite Coxeter groups and their associated simplicial complexes
$\Sigma$.  Details can be found in many places, such as
\cite{\brown,\gb,\humphreys,\tits}.  A \emph{finite reflection group}
on a real inner-product space $V$ is a finite group of orthogonal
transformations of $V$ generated by reflections $s_H$ with respect to
hyperplanes $H$ through the origin.  The set of hyperplanes $H$ such
that $s_H\in W$ is the \emph{reflection arrangement} associated
with~$W$.  Its hyperplane face semigroup $\Sigma$ is called the
\emph{Coxeter complex} of~$W$.  Geometrically, this complex is obtained
by cutting the unit sphere in~$V$ by the hyperplanes~$H$, as in
Section~\ref{sub:cell}.  (As explained there, one might have to first
pass to a quotient of~$V$.)  
It turns out that the Coxeter complex $\Sigma$
is always a simplicial complex.
Furthermore, the action of~$W$ on~$V$ induces an action of
$W$ on~$\Sigma$, and this action is simply-transitive on the chambers.
Thus the set $\C$ of chambers can be identified with~$W$, once a
``fundamental chamber'' $C$ is chosen.

\subsection{Types of simplices}

The number $r$ of vertices of a chamber of $\Sigma$ is called the
\emph{rank} of $\Sigma$ (and of $W$); thus the dimension of $\Sigma$
as a simplicial complex is $r-1$.  It is known that one can color the
vertices of $\Sigma$ with $r$ colors in such a way that vertices
connected by an edge have distinct colors.  The color of a vertex is
also called its \emph{label}, or its \emph{type}, and we denote by $I$
the set of all types.  We can also define $\type(F)$ for any
$F\in\Sigma$; it is the subset of~$I$ consisting of the types of the
vertices of~$F$.  For example, every chamber has type $I$, while the
empty simplex has type $\emptyset$.  The action of $W$ is
type-preserving; moreover, two simplices are in the same $W$-orbit if
and only if they have the same type.  

\subsection{The Coxeter diagram}\label{subs:cd}
Choose a fundamental chamber $C$. It is known that the
reflections $s_i$ in the facets of $C$ generate $W$.
In fact, $W$ has a presentation of the form
$<s_1,\ldots,s_r \mid (s_i s_j)^{m_{ij}}>$
with $m_{ii}=1$ and ${m_{ij}}={m_{ji}}\geq2$.
This data is conveniently encoded in a picture called the Coxeter
diagram of $W$. This diagram is a graph, with vertices and edges,
defined as follows: There are $r$ vertices, one for each generator
$i=1,2,\ldots,r$, and the vertices corresponding to $i$ and $j$ are
connected by an edge if and only if $m_{ij}\geq3$. If $m_{ij}\geq4$
then we simply label the edge with the number $m_{ij}$. 
The figures show the Coxeter diagrams which are of
interest to us, namely the ones of type $A_{n-1}$, $B_n$ and $D_n$.
It is customary to use the generators of $W$, or the vertices of the
Coxeter diagram to label the vertices of
its Coxeter complex $\Sigma$. 
A vertex of the fundamental chamber $C$ is labelled $s_i$ if
it is fixed by all the fundamental reflections except $s_i$.
Since $W$ acts transitively on $\C$ and the action is type-preserving,
this determines the type of all the vertices of $\Sigma$.

\beginpicture
\setcoordinatesystem units <1mm, 1mm>
\setplotarea x from 0 to 80, y from -10 to 20
\put {\small Figure 2.1: Coxeter Diagram of Type $A_{n-1}$} 
  at  40 -5
\circulararc 360 degrees from 11 10 center at 10 10
\put  {$\scriptstyle {s_1}$} at 10 5
\setlinear
\plot 11 10  29 10 /
\circulararc 360 degrees from 31 10 center at 30 10
\put  {$\scriptstyle{s_2}$} at 30 5
\plot 31 10  45 10 /
\put {${\bf \ldots}$} at 50 10
\plot 55 10 69 10 /
\circulararc 360 degrees from 71 10 center at 70 10
\put {$\scriptstyle{s_{n-1}}$} 
 at 70 5
\endpicture

\beginpicture
\setcoordinatesystem units <1mm, 1mm>
\setplotarea x from 0 to 80, y from -10 to 20
\put {\small Figure 2.2: Coxeter Diagram of Type $B_n$} 
  at  50 -5
\circulararc 360 degrees from 11 10 center at 10 10
\put  { $\scriptstyle{s_1}$} at 10 5
\setlinear
\plot 11 10  29 10 /
\circulararc 360 degrees from 31 10 center at 30 10
\put  { $\scriptstyle{s_2}$} at 30 5
\plot 31 10  45 10 /
\put {${\bf \ldots}$} at 50 10
\plot 55 10 69 10 /
\circulararc 360 degrees from 71 10 center at 70 10
\put {$\scriptstyle{s_{n-1}}$} at 70 5
\vanish{
\plot 71 10.5 88 10.5 /
\plot 71 9.5  88 9.5  /
\plot 89 10  87 12 /
\plot 89 10  87  8 /
}
\setlinear
\plot 71 10 89 10 /
\put { $\scriptstyle{4}$}  at 80 15
\circulararc 360 degrees from 91 10 center at 90 10
\put {$\scriptstyle{s_n=t}$} 
 at 90 5
\endpicture

\beginpicture
\setcoordinatesystem units <1mm, 1mm>
\setplotarea x from 0 to 80, y from -10 to 25
\put {\small Figure 2.3: Coxeter Diagram of Type $D_n$} 
  at  50 -10
\circulararc 360 degrees from 11 10 center at 10 10
\put  { $\scriptstyle{s_1}$} at 10 5
\setlinear
\plot 11 10  29 10 /
\circulararc 360 degrees from 31 10 center at 30 10
\put  { $\scriptstyle{s_2}$} at 30 5
\plot 31 10  45 10 /
\put {${\bf \ldots}$} at 50 10
\plot 55 10 69 10 /
\circulararc 360 degrees from 71 10 center at 70 10
\put { $\scriptstyle{s_{n-2}}$} 
 at 70 5
\circulararc 360 degrees from 86 20 center at 85 20
\circulararc 360 degrees from 86 0 center at 85 0
\plot 71 10.5  84 19.5 /
\plot 71 9.5  84 0.5 /
\put { $\scriptstyle{u}$}  at 87 16
\put { $\scriptstyle{v}$}  at 87 -4
\endpicture

\subsection{The Coxeter group of type $A_{n-1}$} \label{subs:a}

The Coxeter group $W=S_n$ acts on $\R^n$ by permuting the
coordinates. 
The arrangement in this case is the \emph{braid
arrangement} in $\R^n$. It is discussed in detail in
\cite{\bid,\bhr,\bbd,\bd}. 
It consists of the $\binom{n}{2}$
hyperplanes $H_{ij}$ defined by $x_i=x_j$, where $1\le i<j\le n$.
Each chamber is determined by an ordering of the coordinates, so it
corresponds to a permutation.  
The faces of a chamber are obtained by
changing to equalities some of the inequalities defining that chamber.

We fix  
$x_1<x_2<\ldots<x_n$ 
to be the fundamental chamber $C$.
The supports of the facets of $C$ are hyperplanes of the form
$x_i=x_{i+1}$, where $1\le i \le n-1$.
The reflection in the hyperplane $x_i=x_{i+1}$
correponds to the generator $s_i$ that
interchanges the coordinates $x_i$
and $x_{i+1}$.
The chamber $C$
has $n-1$ vertices, namely

$s_1 : x_1<x_2=\ldots=x_n$, 

$s_2 : x_1=x_2<x_3=\ldots=x_n$,$\ldots$, 

$s_{n-1} : x_1=\ldots=x_{n-1}<x_n$.

\noindent The labels $s_1$,$s_2$,$\ldots$,$s_{n-1}$ are assigned by the rule
mentioned in subsection~\ref{subs:cd}. Applying the action of $W$ we see,
for example, that $x_{\pi(1)}<x_{\pi(2)}=\ldots=x_{\pi(n)}$ gives all
vertices of type $s_1$ as $\pi$ varies over all permutations of~$[n]=\{1,\dots,n\}$. 
\vanish{
When $n=4$, for example, one of the 24
chambers is the region defined by $x_2>x_3>x_1>x_4$, corresponding to
the permutation 2314.  
For
example, the chamber $x_2>x_3>x_1>x_4$ has a face given by
$x_2>x_3>x_1=x_4$, which is also a face of the chamber
$x_2>x_3>x_4>x_1$.
}

We encode the system of equalities and inequalities
defining a face~$F$ by an ordered partition $(B_1,\dots,B_k)$ of
$[n]$.  Here $B_1,\dots,B_k$ are disjoint nonempty sets
whose union is $[n]$, and their order counts.  
\vanish{
For example, the face
$x_2>x_3>x_1=x_4$ corresponds to the 3-block ordered partition
$(\{2\},\{3\},\{1,4\})$, and the face $x_2>x_1=x_3=x_4$ corresponds to
the 2-block ordered partition $(\{2\}, \{1,3,4\})$.
}
Thus the simplices of $\Sigma$ are ordered partitions $B=(B_1,\dots,B_l)$ 
of the set $[n]$. 
For example, the vertices of type $s_1$ are ordered two block partitions
such that the first block is a singleton.

The product in $\Sigma$ is also easy to
describe. We multiply two ordered
partitions by taking intersections and ordering them
lexicographically; more precisely, if $B=(B_1,\dots,B_l)$ and
$C=(C_1,\dots,C_m)$, then
\[
BC=(B_1\cap C_1,\dots,B_1\cap C_m,\dots, B_l\cap C_1,\dots,B_l\cap
C_m)\sphat\,,
\]    
where the hat means ``delete empty intersections''. The 1-block
ordered partition is the identity.

Note that $B$ is a face of $C$ if and only if $C$ consists of
an ordered partition of $B_1$ followed by an ordered partition
of~$B_2$, and so on, that is, if and only if $C$ is a refinement of
$B$. The chambers are the ordered partitions into
singletons, so they correspond to the permutations of~$[n]$ or a deck
of $n$ cards.

\vanish{
It is useful to have a second description of $\B$.  Ordered partitions
$(B_1,\dots,B_l)$ of~$[n]$ are in 1--1 correspondence with chains of
subsets $\emptyset=E_0<E_1<\cdots<E_l=[n]$, the correspondence being
given by $B_i=E_i-E_{i-1}$.  So we may identify $\B$ with the set of
such chains.  
The
vertices are the proper nonempty subsets $X\subset[n]=\{1,\dots,n\}$,
and the simplices are the chains of such subsets.  The $S_n$-action is
induced by the action of $S_n$ on~$[n]$. The
chambers of~$\Sigma$ correspond to permutations $w$ of $[n]$, with $w$
corresponding to the chamber
\[
\{w(1)\}<\{w(1),w(2)\}<\cdots<\{w(1),w(2),\dots,w(n-1)\}.
\]
This is the same as the identification of $\C$ with $S_n$ that results
from choosing
\[
\{1\}<\{1,2\}<\cdots<\{1,2,\dots,n-1\}
\]
as fundamental chamber.
The product is then described as follows:  Given a chain
$E$ as above and a second chain $F\colon\emptyset=F_0<F_1<\cdots<
F_m=[n]$, their product $EF$ is obtained by using $F$ to refine~$E$.
More precisely, consider the sets $G_{ij}=(E_{i-1}\cup F_j)\cap E_i =
E_{i-1}\cup (F_j\cap E_i)$.  For each $i=1,2,\dots,l$ we have
\[
E_{i-1} = G_{i0}\subseteq G_{i1}\subseteq\cdots\subseteq G_{im}=E_i.
\]
Deleting repetitions gives a chain from $E_{i-1}$ to~$E_i$, and
combining these for all~$i$ gives the desired refinement $EF$ of~$E$.

This construction is used in one of the standard proofs of the
Jordan--H\"{o}lder theorem.
}

\subsection{The Coxeter group of type $B_n$} \label{subs:b}
 
The group of signed permutations $W=S_n \ltimes \Z_2^n$ acts on $\R^n$
with the subgroup $S_n$ permuting the
coordinates and the subgroup $\Z_2^n$ flipping the signs of the
coordinates. 
The reflection arrangement in this case consists of the hyperplanes
defined by $x_i= \pm x_j$ and $x_i = 0$, where $1\le i<j\le n$.
A chamber is given by an ordering of the coordinates, their negatives
and zero. For example, for $n = 3$, 
$x_2<-x_3<x_1<0<-x_1<x_3<-x_2$ specifies a chamber.
Note that the inequalities that appear on the left of $0$ completely
determine a chamber (and the same is true for any face since it is obtained by changing to equalities some of the inequalities defining a chamber).
Thus we see that a chamber corresponds to a signed permutation.

We fix
$x_1<x_2<\ldots<x_n<0$ 
to be the fundamental chamber $C$. 
The supports of the facets of $C$ are hyperplanes of the form
$x_i=x_{i+1}$, where $1\le i \le n-1$ and $x_n = 0$.
The generators $s_i$ interchange the coordinates $x_i$
and $x_{i+1}$ and the generator $t$ flips the sign of $x_n$.
The $n$ vertices of $C$ along with their labels are as follows.

$s_1:x_1<x_2=\ldots=x_n=0$, 

$s_2:x_1=x_2<x_3=\ldots=x_n=0$,$\ldots$, 

$s_{n-1}:x_1=x_2=\ldots=x_{n-1}<x_n=0$,

$t:x_1=x_2\ldots=x_n<0$.

\noindent Applying the group action, we see, for example, that the
vertices of type $t$ have the form
$\varepsilon_1 x_1=\varepsilon_2 x_2\ldots=\varepsilon_n
x_n<0$
where $\varepsilon_i \in \{\pm 1 \}$.

It is convenient to describe a face~$F$ by an anti-symmetric ordered partition
$(B_1,\dots,B_k,Z,\overline{B}_k,\ldots,\overline{B}_1)$ of
$[n,\overline{n}]=\{1,\dots,n,\overline{1},\ldots,\overline{n}\}$. 
We will call this a {\it partition of type $B_n$}.
For example, for $n = 3$, the face
$x_2<-x_3<x_1=0=-x_1<x_3<-x_2$ is written
$(\{2\},\{\overline{3}\},\{1,\overline{1}\},\{3\},\{\overline{2}\})$.
The set $Z$ satisfies $Z = \overline{Z}$ and is allowed to be
empty. We call it the zero block. It is the only set that is allowed
to contain both a number and its negative. 
Also, the sets $B_1,\dots,B_k$ are necessarily non-empty.
We also define a \emph{weak partition of type $B_n$} to be a partition as
above but where
the sets $B_1,\dots,B_k$ are allowed to be empty.
Note that for a face $F$, the zero block of its partition is empty if
and only if the type of $F$ contains the letter 
$t$. 
We split the description for a face $F$
into four cases, depending on whether the type of $F$ contains 
\begin{itemize}
\item[(i)] neither $s_{n-1}$ nor $t$: 

An ordered partition
$(B_1,\dots,B_k,Z,\overline{B}_k,\ldots,\overline{B}_1)$ of
$[n,\overline{n}]$ of type $B_n$ with the added
restriction that the zero block $Z$ has at least $4$ elements.

\item[(ii)] $t$ but not $s_{n-1}$:

An ordered partition
$(B_1,\dots,B_k,Z,\overline{B}_k,\ldots,\overline{B}_1)$ of
$[n,\overline{n}]$, with the restriction that the zero block is empty
and $B_k$ has at least $2$ elements. 

\item[(iii)] $s_{n-1}$ but not $t$:

An ordered partition
$(B_1,\dots,B_k,Z,\overline{B}_k,\ldots,\overline{B}_1)$ of
$[n,\overline{n}]$, with the restriction that the zero block has exactly
$2$ elements. 

\item[(iv)] both $s_{n-1}$ and $t$:

An ordered partition
$(B_1,\dots,B_k,Z,\overline{B}_k,\ldots,\overline{B}_1)$ of
$[n,\overline{n}]$, with the restriction that the zero block is empty
and $B_k$ has at exactly $1$ element. 

\end{itemize}
The product in $\Sigma$ is defined exactly as for the $A_{n-1}$ case,
where we refine the first partition by the second. 
Note that $B$ is a face of $C$ if and only if $C$ is a refinement of
$B$.

As mentioned earlier, the second half of the partition following $Z$
contains no new information. Hence we may describe a $k-1$ dimensional
face by a $(k+1)$ block partition 
$(B_1,\dots,B_k,Z)$ of $[n]$. We think of $B_1,\dots,B_k$ as signed sets
and $Z$ (possibly empty) as an unsigned set.
The chambers are then $(n+1)$ block
partitions into $n$ singletons and an empty zero block, so they correspond
to signed permutations of~$[n]$ or a deck of $n$ signed cards.
We will call such a deck as a \emph{deck of type $B_n$}.
However, the product and the face relation is not so natural now. So we
will mainly use the first description.

\subsection{The Coxeter group of type $D_n$} \label{subs:d}

The group of even signed permutations $W=S_n \ltimes G$ is an index 2
subgroup of the group of signed permutations $S_n \ltimes \Z_2^n$.
Here $G$ is the index 2 subgroup of $\Z_2^n$ consisting of $n$-tuples
that have an even number of negative signs.
The group $W$ acts on $\R^n$
with the subgroup $S_n$ permuting the
coordinates and the subgroup $G$ flipping the signs of the
coordinates. The
reflection arrangement in this case consists of the hyperplanes
defined by $x_i= \pm x_j$, where $1\le i<j\le n$.
It is obtained from the reflection arrangement of type $B_n$ by
deleting the coordinate hyperplanes $x_i=0$ for $1\le i \le n$.
On the level of chambers, this has the effect of merging together pairs
of adjacent chambers of the Coxeter complex of type $B_n$; 
see Figure~\ref{maps}.
For example, the chamber of the Coxeter complex of type $D_n$,
$x_1<x_2<\ldots<x_{n-1}<\pm x_n<-x_{n-1}<\ldots<-x_1$
(which we fix as our fundamental chamber $C$) is the union of the two
chambers of the Coxeter complex of type $B_n$, namely, 
$x_1<x_2<\ldots<x_n<0$ and $x_1<x_2<\ldots<-x_n<0$.
The supports of the facets of $C$ are hyperplanes of the form
$x_i=x_{i+1}$, where $1\le i \le n-1$ and $x_{n-1} = - x_n$.
The generators $s_i$ interchange the coordinates $x_i$
and $x_{i+1}$, the generator $u$ interchanges the coordinates $x_{n-1}$
and $x_{n}$ and the generator $v$ interchanges $x_{n-1}$ and $x_n$ and
flips the signs of both.
The $n$ vertices of $C$ along with their labels are as follows.

$s_1:x_1<x_2=\ldots=x_n=0$, 

$s_2:x_1=x_2<x_3=\ldots=x_n=0$,$\ldots$, 

$s_{n-2}:x_1=x_2=\ldots=x_{n-2}<x_{n-1}=x_n=0$,

$u:x_1=x_2\ldots=x_{n-1}=-x_n<x_n=-x_{n-1}=\ldots=-x_{2}=-x_1$,

$v:x_1=x_2\ldots=x_{n-1}=x_n<-x_n=-x_{n-1}=\ldots=-x_{2}=-x_1$.

\noindent Applying the group action, we see, for example, that the
vertices of type $u$ and $v$ both have the form
$\varepsilon_1 x_1=\varepsilon_2 x_2\ldots=\varepsilon_n
x_n<-\varepsilon_n x_n=\ldots=-\varepsilon_2
x_{2}=-\varepsilon_1 x_1$, 
where $\varepsilon_i \in \{\pm 1 \}$
and they are distinguished by the parity of the product 
$\varepsilon_1 \varepsilon_2 \ldots \varepsilon_n$.
Also observe that the edge of the fundamental chamber $C$ of type
$\{u,v\}$ is given by 
$x_1=x_2\ldots=x_{n-1}<\pm x_n<-x_{n-1}=\ldots=-x_{2}=-x_1$.

As in the $B_n$ case, we 
describe a face~$F$ by an anti-symmetric ordered partition
$(B_1,\dots,B_k,C,\overline{B}_k,\ldots,\overline{B}_1)$ of
$[n,\overline{n}]$.
We repeat that the sets $B_1,\dots,B_k$ are
non-empty and cannot contain both a number and its negative. 
We call the block $C=\overline{C}$ in the middle, the central
block rather than the zero block. 
It always has an even number of elements.
Furthermore, we impose the relation
$(B_1,\dots,B_{k-1},C,\overline{B}_{k-1},\ldots,\overline{B}_1) =
(B_1,\dots,B_{k-1},B_{k},C^{\prime},\overline{B}_{k},\overline{B}_{k-1},\ldots,\overline{B}_1)$ 
where $C^{\prime}$ is empty, $C$ has exactly two elements, namely, a number and
its negative, and $C=B_{k} \cup C^{\prime} \cup \overline{B}_{k}$. 
We split the description for a face $F$
into three cases, depending on whether the type of $F$ contains 
\begin{itemize}
\item[(i)] neither $u$ nor $v$: 

An ordered partition
$(B_1,\dots,B_k,C,\overline{B}_k,\ldots,\overline{B}_1)$ of
$[n,\overline{n}]$, with the
restriction that the central block $C$ (in this case, we may think of
it as the zero block) has at least $4$ elements.

\item[(ii)] exactly one of $u$ and $v$:

An ordered partition
$(B_1,\dots,B_k,C,\overline{B}_k,\ldots,\overline{B}_1)$ of
$[n,\overline{n}]$, with the restriction that $C$ is empty and $B_k$
has at least $2$ elements. 

\item[(iii)] both $u$ and $v$:

An ordered partition
$(B_1,\dots,B_{k-1},C,\overline{B}_{k-1},\ldots,\overline{B}_1)$ of
$[n,\overline{n}]$, with the restriction that $C$ has exactly $2$
elements, namely, a number and its negative.  

An ordered partition
$(B_1,\dots,B_{k-1},B_{k},C^{\prime},\overline{B}_{k},\overline{B}_{k-1},\ldots,\overline{B}_1)$ 
where $C^{\prime}$ is empty and $B_{k}$ has exactly $1$ element.

\end{itemize}
Here the number $k$ is always the rank of the
face. 
The first three descriptions are obtained directly from the equalities
and inequalities that define a face. The reason why we call $C$ the
central block rather than the zero block should be clear now.
The motivation for the second description in case $(iii)$ is as follows. 
To give an example, the product of
$(\{2,\overline{3}\},\{1,\overline{5},\overline{4}\},\{4,5,\overline{1}\},\{3,\overline{2}\})$
with
$(\{4,1,\overline{5}\},\{2,\overline{3},3,\overline{2}\},\{5,\overline{1},\overline{4}\})$,
must be
$(\{2,\overline{3}\},\{1,\overline{5}\},\{\overline{4},4\},\{5,\overline{1}\},\{3,\overline{2}\})$.
But if we refine the first partition by the second then
we get
$(\{2,\overline{3}\},\{1,\overline{5}\},\{\overline{4}\},\{4\},\{5,\overline{1}\},\{3,\overline{2}\})$.
This partition has two singletons in the middle and does not fit any
of the first three descriptions given above. So we allow ourselves to merge the
two singletons in the middle.
This identification allows us to multiply two partitions $B$ and $C$
in the same way as before; that 
is, we refine $B$ using $C$. 

As in the $B_n$ case, if we throw off the (redundant) second half of
the partition then we see that chambers correspond to an \emph{almost}
signed deck of cards or a \emph{deck of type $D_n$}. It is a deck in
which every card, except the bottomost, is signed.
We also define a \emph{weak partition of type $D_n$} to be a partition 
of type $D_n$ where
the sets $B_1,\dots,B_k$ are allowed to be empty.

\subsection*{Acknowledgments} It is a pleasure to thank my advisor
Ken Brown for introducing me to this area of mathematics and sharing
his insights with me. I would also like to thank Marcelo Aguiar for
giving valuable suggestions.



\providecommand{\bysame}{\leavevmode\hbox to3em{\hrulefill}\thinspace}
\providecommand{\MR}{\relax\ifhmode\unskip\space\fi MR }
\providecommand{\MRhref}[2]{%
  \href{http://www.ams.org/mathscinet-getitem?mr=#1}{#2}
}
\providecommand{\href}[2]{#2}

\end{document}